\documentclass[11pt,letterpaper]{amsart}
\usepackage[T1]{fontenc}
\usepackage{lmodern}
\usepackage{amsmath,amssymb,amsthm,mathtools,mathrsfs}
\usepackage{enumitem}
\usepackage{microtype}
\usepackage{tikz}
\usetikzlibrary{arrows.meta}

\usepackage{xurl}
\usepackage[bookmarksnumbered,hidelinks]{hyperref}
\hypersetup{pdftitle={Two-scale band geometry and ballistic transport for limit-periodic Schrodinger operators},
  pdfauthor={Moises Gomez Solis},
  pdfsubject={Limit-periodic Schrodinger operators and ballistic transport},
  pdfkeywords={limit-periodic, Schrodinger operators, ballistic transport, Floquet theory}}
\numberwithin{equation}{section}
\mathtoolsset{showmanualtags}
\newtheorem{theorem}{Theorem}[section]
\newtheorem{lemma}[theorem]{Lemma}
\newtheorem{proposition}[theorem]{Proposition}
\newtheorem{corollary}[theorem]{Corollary}
\theoremstyle{definition}
\newtheorem{definition}[theorem]{Definition}
\newtheorem{example}[theorem]{Example}
\theoremstyle{definition}
\newtheorem{remark}[theorem]{Remark}
\newtheorem*{localremark}{Remark}
\newcommand{\R}{\mathbb R}

\newcommand{\Z}{\mathbb Z}
\newcommand{\N}{\mathbb N}
\newcommand{\T}{\mathbb T}
\newcommand{\al}{\alpha}
\newcommand{\e}{\varepsilon}
\newcommand{\Bz}{\mathcal B}
\newcommand{\EC}{\mathrm{EC}}
\newcommand{\ECY}{\mathrm{EC}_{\mathrm Y}}
\newcommand{\Hk}{H(k)}
\newcommand{\Hz}{H_0}
\newcommand{\dd}{\mathrm{d}}
\newcommand{\dkn}[1]{\frac{\dd k}{|\T^*_{#1}|}}
\newcommand{\dk}{\frac{\dd k}{|\T^*|}}
\newcommand{\Dom}{\operatorname{Dom}}
\newcommand{\mes}{\operatorname{mes}}
\newcommand{\eqdef}{\overset{\mathrm{def}}{=\joinrel=}}
\DeclareMathOperator{\dist}{dist}
\title[Two-scale band geometry and ballistic transport]{Two-scale band geometry and ballistic transport for limit-periodic Schr\"odinger operators}
\author[M. Gomez Solis]{Moises Gomez Solis}
\address{Department of Mathematics, Louisiana State University,
303 Lockett Hall, 175 Field House Drive,
Baton Rouge, LA 70803-4918, USA}
\email{mgome29@lsu.edu}
\subjclass[2020]{Primary 34L40; Secondary 35J10, 81Q10}
\keywords{Limit-periodic Schr\"odinger operators, ballistic transport, Floquet theory,
asymptotic velocity, periodic approximation}
\date{}

\DeclareMathOperator{\tr}{tr}
\DeclareMathOperator{\Ran}{Ran}
\begin{document}
\begin{abstract}
We prove strong ballistic transport for one-dimensional Schr\"odinger
operators with bounded limit-periodic potentials. If the nested periodic
approximants satisfy
$\|V-V_n\|_\infty=o(p_{n+1}^{-2}e^{-\eta p_n})$, with
$\eta>7\sqrt{2R}$ and $\|V\|_\infty\le R$, then the rescaled Heisenberg
position has a nonzero strong limit on a weighted $H^2$ class containing
the Schwartz space. In contrast to earlier continuum results requiring
exponential accuracy in $p_{n+1}$, the exponential weight here involves
only $p_n$; the ratios $p_{n+1}/p_n$ may be unbounded. Folding identifies
each ordered large-cell band, away from symmetry momenta, with a fixed
translated small-cell band. This confines near-collisions of all bands to
a common neighborhood controlled at the smaller period. A period-uniform
quadratic-form estimate, ordered-gap summation, and a graph-norm comparison
linear in the potential error then give the limiting velocity. Moments of
orders $0<r\le4$ have two-sided ballistic growth as a corollary. An explicit
example separates the two approximation conditions.
\end{abstract}
\maketitle

\section{Introduction}\label{sec:intro}
Periodic quantum dynamics admit a nonzero asymptotic velocity under the
hypotheses of Asch and Knauf \cite{AsKn}. For limit-periodic operators,
quantitative periodic approximation transfers this behavior to the limit;
see \cite{DYL,Fil2017} for Jacobi and block Jacobi matrices and \cite{Young}
for continuum Schr\"odinger operators. The two-scale argument of Han and
Gomez Solis \cite{HGS} makes the smaller of two consecutive periods control
the costly spectral estimate. The present paper develops this mechanism in
the continuum, with an explicit sufficient approximation threshold.
For the notion of strong ballistic transport used here, see
\cite[Definition 2.1(1)]{DMY26}.

The starting point is geometric. Representing a $p$-periodic potential in a
cell of length $P=\ell p$ does not change the whole-line operator or its
spectrum. It changes the fiber representation. Translation by $p$ splits
the $P$-cell space into $\ell$ character sectors, each carrying a translated
copy of a small-cell fiber. The resulting eigenvalues must still be ordered
as the spectrum of a genuine $P$-periodic scalar operator. Away from the
large-cell symmetry momenta, simplicity and continuity keep the small-cell
label of each ordered band fixed. This allows us to combine large-cell band
geometry with derivative estimates at the smaller period.

There are two separate comparisons. First, folding is an exact unitary
decomposition of the old operator. Second, a bounded perturbation compares
its ordered eigenvalues with those of the new operator in the same large
cell. The first comparison retains small-period slopes; the second moves
each ordered eigenvalue by at most the potential error. Together they yield
band-width bounds and a common localization of near-collisions, uniformly
in the band index. The near-collision estimate enters the periodic averaging
argument and, through a linear dynamical comparison, the limit theorem.

\subsection{The approximation class and main result}
Let $V:\R\to\R$ be bounded and measurable. On $L^2(\R)$ we consider
the self-adjoint operator
\begin{equation}\label{eq:H}
 H\eqdef-\frac{\dd^2}{\dd x^2}+V,\qquad \Dom(H)=H^2(\R).
\end{equation}
All potential norms are essential supremum norms. The evolution
$e^{-itH}\psi$ solves the Schr\"odinger equation for $\psi\in H^2(\R)$.
We write $X$ for multiplication by $x$, with
\begin{equation}\label{def. position operator}
 \Dom(X)=\{\psi\in L^2(\R):x\psi\in L^2(\R)\},
\end{equation}
and write
\begin{equation}\label{def. heisenberg evolution}
 X_H(t)\eqdef e^{itH}Xe^{-itH}
\end{equation}
on vectors for which $e^{-itH}\psi\in\Dom(X)$.

\begin{definition}[Two-scale exponential class]\label{def:class}
Let $\eta,\alpha>0$. We write $H\in\EC(\eta,\alpha)$ if there are
continuous real $p_n$-periodic potentials $V_n$, with $p_n>0$ and
$p_{n+1}/p_n\in\N\setminus\{1\}$, such that
\begin{equation}\label{eq. exponential class}
 \lim_{n\to\infty}p_{n+1}^{\alpha}e^{\eta p_n}\|V-V_n\|_\infty=0.
\end{equation}
The displayed periods need not be minimal.
\end{definition}
The parameter $\alpha$ is retained to compare the main case $\alpha=2$
with the power $4$ obtained from a square-root perturbation estimate.
For $s>0$, set
\begin{equation}\label{eq:Ds}
 \mathcal D_s\eqdef
 \{\psi\in H^2(\R):|x|^s\psi\in L^2(\R),\ |x|^s\psi''\in L^2(\R)\}.
\end{equation}
This space contains the Schwartz space and, for $s\ge2$, is contained in
$H^2(\R)\cap\Dom(X^2)$.

\begin{theorem}[Strong ballistic transport]\label{thm:main}
Let $R>0$ and $\eta>7\sqrt{2R}$. Suppose that $\|V\|_\infty\le R$ and
$H\in\EC(\eta,2)$. For every $s>2$ and every
$\psi\in\mathcal D_s\setminus\{0\}$, the limit
\begin{equation}\label{eq:mainlimit}
 Q\psi\eqdef\lim_{t\to\infty}t^{-1}X_H(t)\psi
\end{equation}
exists in $L^2(\R)$ and is nonzero.
\end{theorem}
Clipping the approximants to $[-R,R]$ preserves the approximation rate;
see Section~\ref{sec:fold}. Thus one may take $R=\|V\|_\infty$ when
$V\ne0$. The free case is already periodic and can also be treated with
any positive bound $R$.

\subsection{Comparison with the preceding continuum theorem}
The role of the approximation rate already appears at the spectral level.
In the continuum, purely singular continuous spectrum on a zero-measure
Cantor set is generic in a natural limit-periodic space \cite{DFLthin}.
The Pastur--Tkachenko condition gives purely absolutely continuous
spectrum \cite{PT84,PT88}, and this spectrum is homogeneous
\cite[Theorem 1.1]{Fil2015}. For transport, the comparison is with
Young's continuum theorem.

To distinguish the two conventions, write $H\in\ECY(\eta')$ when the
same type of nested periodic approximation satisfies
\begin{equation}\label{eq:young-class}
 e^{\eta'p_{n+1}}\|V-V_n\|_\infty\longrightarrow0.
\end{equation}
Young's Theorem~1.2 \cite{Young} proves the nonzero strong limit on the
same spaces $\mathcal D_s$, $s>2$, for $\eta'$ sufficiently large in
terms of $R$. His threshold is not given as an explicit numerical function
of $R$ in the theorem.

For the question of optimal approximation rates, see \cite{DF2,HGS}.
In dimension two, Karpeshina, Lee, Shterenberg, and Stolz
\cite{KLSS} prove time-averaged ballistic lower bounds for suitable
initial states and certain classes of limit-periodic or quasi-periodic
potentials.

\begin{remark}[Comparison of approximation conditions]\label{rem:young-comparison}
For $0<\eta<2\eta'$, the inequality $p_{n+1}\ge2p_n$ gives
\begin{equation}\label{eq:class-inclusion}
 \begin{aligned}
 p_{n+1}^2e^{\eta p_n}\|V-V_n\|_\infty
 &\le p_{n+1}^2e^{-(\eta'-\eta/2)p_{n+1}}
       \bigl(e^{\eta'p_{n+1}}\|V-V_n\|_\infty\bigr)
 \longrightarrow0.
 \end{aligned}
\end{equation}
Hence $\ECY(\eta')\subseteq\EC(\eta,2)$ in this range. In particular,
$\eta'>\tfrac72\sqrt{2R}$ is sufficient to apply
Theorem~\ref{thm:main} through this inclusion. Conversely, if a witnessing
sequence satisfies $p_{n+1}/p_n\le L$, its two-scale condition implies
\eqref{eq:young-class} for every $\eta'<\eta/L$. For unbounded ratios
there is a genuine distinction, as the following example shows.
The constants in the two discriminant estimates are different; the
comparison here is of the stated hypotheses, not a fixed numerical ratio
between their proofs' exponential constants.
\end{remark}

\begin{example}[A potential outside every next-period exponential class]
\label{ex:strictclass}
Fix $\lambda>0$, put $p_1=4$, and define
\[
 \begin{gathered}
 p_{k+1}=p_k\left\lceil\frac{e^{p_k}}{p_k}\right\rceil,
 \qquad a_1=\lambda/2,\quad a_{k+1}=\lambda e^{-p_k^2},\\
 V(x)=\sum_{k\ge1}a_k\cos(2\pi x/p_k).
 \end{gathered}
\]
Then $p_k\mid p_{k+1}$, $p_{k+1}\ge2p_k$, and
$e^{p_k}\le p_{k+1}\le2e^{p_k}$. The series converges uniformly,
$\|V\|_\infty<\lambda$, and its partial sums $V_n$ satisfy
\[
 \|V-V_n\|_\infty\le2\lambda e^{-p_n^2},\qquad
 p_{n+1}^2e^{\eta p_n}\|V-V_n\|_\infty
 \le8\lambda e^{-p_n^2+(\eta+2)p_n}\longrightarrow0
\]
for every $\eta>0$. The tail bound follows because successive amplitudes
have ratio at most $1/2$.

This potential belongs to no $\ECY(\eta')$, even with a different
approximating sequence. Indeed, for a uniformly almost-periodic function
$f$, its mean coefficient at frequency $\xi$ is
(see \cite[Chapter~I, \S3]{Besic})
\[
 \widehat f_{\rm B}(\xi)
 =\lim_{T\to\infty}\frac1{2T}\int_{-T}^{T}f(x)e^{-2\pi i\xi x}\,\dd x,
 \qquad |\widehat f_{\rm B}(\xi)|\le\|f\|_\infty.
\]
Uniform convergence permits computation term by term, giving
$\widehat V_{\rm B}(1/p_{k+1})=a_{k+1}/2$.
A continuous $q$-periodic function has zero coefficient at this frequency
when $q<p_{k+1}$: translation of the defining mean by $q$ multiplies it
by $e^{2\pi iq/p_{k+1}}\ne1$.
For any competing nested periods $q_m\to\infty$, choose the last $m=m(k)$
with $q_m<p_{k+1}$, so $q_{m+1}\ge p_{k+1}$ and $m(k)\to\infty$.
Every corresponding $q_m$-periodic approximant $V_m^\sharp$ obeys
\[
 e^{\eta'q_{m+1}}\|V-V_m^\sharp\|_\infty
 \ge\frac\lambda2\exp(\eta'e^{p_k}-p_k^2)\longrightarrow\infty.
\]
The same nonzero coefficients at frequencies tending to zero exclude any
fixed period for $V$. The example is therefore strictly limit-periodic.
\end{example}

Write $\ECY(\infty)=\bigcap_{\eta'>0}\ECY(\eta')$, with the
approximating sequence allowed to depend on $\eta'$.
Requiring one sequence of continuous periodic approximants to satisfy
\eqref{eq:young-class} for every $\eta'>0$ is the uniform-norm form of
the Pastur--Tkachenko condition \cite{Fil2015,PT84,PT88}; it implies
membership in $\ECY(\infty)$. The inclusion
\eqref{eq:class-inclusion} immediately recovers the following consequence.
\begin{corollary}\label{cor:young-infty}
If $H\in\ECY(\infty)$ and $V$ is bounded, the conclusion of
Theorem~\ref{thm:main} holds for every nonzero $\psi\in\mathcal D_s$, $s>2$.
\end{corollary}
This consequence is also contained in \cite[Theorem 1.2]{Young}.

\subsection{Moment growth}
For a nonzero state, define
\begin{equation}\label{def. p-moments}
 M_{r,\psi}(t)\eqdef\int_\R |x|^r|e^{-itH}\psi(x)|^2\,\dd x,
 \qquad r>0,
\end{equation}
and, when these moments are finite for all sufficiently large $t$, set
\begin{equation}\label{def. transport exponents}
 \begin{aligned}
 \beta^+_\psi(r)&\eqdef\limsup_{t\to\infty}
       \frac{\log M_{r,\psi}(t)}{r\log t},\\
 \beta^-_\psi(r)&\eqdef\liminf_{t\to\infty}
       \frac{\log M_{r,\psi}(t)}{r\log t}.
 \end{aligned}
\end{equation}
\begin{corollary}[Ballistic moments]\label{cor:moments}
Under the hypotheses of Theorem~\ref{thm:main}, for every
$\psi\in\mathcal D_s\setminus\{0\}$, $s>2$, and $0<r\le4$, there are
$a_{r,\psi},b_{r,\psi}>0$ and $T_{r,\psi}<\infty$ such that
\begin{equation}\label{eq:ballisticmoments}
 a_{r,\psi}t^r\le M_{r,\psi}(t)\le b_{r,\psi}t^r,
 \qquad t\ge T_{r,\psi}.
\end{equation}
In particular, $\beta^-_\psi(r)=\beta^+_\psi(r)=1$ in this range.
\end{corollary}
This is an interpolation consequence of the velocity limit and the
fourth-moment estimate of Radin--Simon \cite{RadSim}, in the precise
bounded-potential form proved by Young \cite[Theorem 5.2]{Young}.
The same consequence follows from Young's transport theorem on his
approximation class. A proof is given in Section~\ref{sec:main}; the
endpoint restriction for the present state class is discussed in
Section~\ref{sec:remarks}.

\subsection{The geometric estimate and the proof}
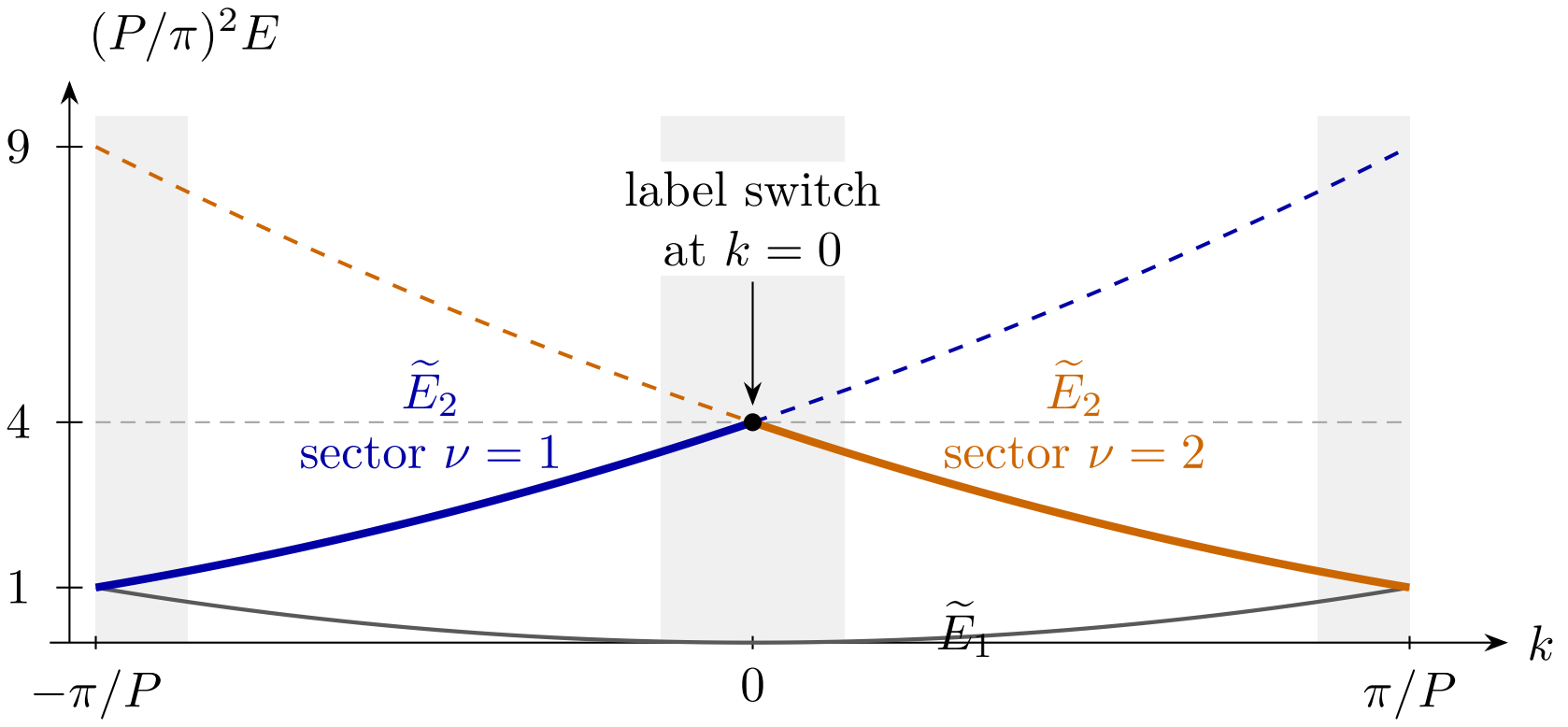
\begin{figure}[t]
\centering
\begin{tikzpicture}[x=4.65cm,y=.39cm,>=Stealth,font=\small]
 \fill[gray!12] (-1,0) rectangle (-.86,9.55);
 \fill[gray!12] (-.14,0) rectangle (.14,9.55);
 \fill[gray!12] (.86,0) rectangle (1,9.55);
 \draw[->] (-1.07,0)--(1.15,0) node[right] {$k$};
 \draw[->] (-1.04,0)--(-1.04,10.2) node[above right] {$ (P/\pi)^2E$};
 \foreach \x/\lab in {-1/{-\pi/P},0/{0},1/{\pi/P}} {
   \draw (\x,-.12)--(\x,.12);
   \node[below] at (\x,-.14) {$\lab$};
 }
 \foreach \y in {1,4,9} {
   \draw (-1.06,\y)--(-1.02,\y);
   \node[left] at (-1.07,\y) {$\y$};
 }
 \draw[gray!70,densely dashed] (-1,4)--(1,4);
 \draw[black!65,line width=.8pt,domain=-1:1,samples=70] plot (\x,{\x*\x});
 \draw[blue!65!black,line width=.75pt,dashed,domain=0:1,samples=50] plot (\x,{(\x+2)^2});
 \draw[orange!80!black,line width=.75pt,dashed,domain=-1:0,samples=50] plot (\x,{(\x-2)^2});
 \draw[blue!65!black,line width=1.6pt,domain=-1:0,samples=50] plot (\x,{(\x+2)^2});
 \draw[orange!80!black,line width=1.6pt,domain=0:1,samples=50] plot (\x,{(\x-2)^2});
 \fill (0,4) circle (1.8pt);
 \node[anchor=west] at (.25,.25) {$\widetilde E_1$};
 \node[text=blue!65!black,align=center] at (-.49,4.15)
   {$\widetilde E_2$\smallskip\\sector $\nu=1$};
 \node[text=orange!80!black,align=center] at (.49,4.15)
   {$\widetilde E_2$\smallskip\\sector $\nu=2$};
 \node[fill=white,inner sep=2pt,align=center] at (0,7.7)
   {label switch\\at $k=0$};
 \draw[->] (0,6.55)--(0,4.3);
\end{tikzpicture}
\caption{Threefold enlargement of the free operator, $P=3p$.
The solid colored curve is the second ordered band
$\widetilde E_2(k)=(2\pi/P-|k|)^2$; the dashed curves continue its
individual small-cell branches. Its sector is constant on each open
half-zone and changes at $k=0$. The gray strips depict neighborhoods of
$0,\pm\pi/P$, where the two-scale estimate permits near-collisions
after a small perturbation. Their widths are schematic; the curves are exact.}
\label{fig:folding}
\end{figure}

We now make the folding mechanism quantitative. Fix $n$ and put
$p=p_n$ and $P=p_{n+1}=\ell p$. With the notation of
Section~\ref{sec:fold}, each ordered large-cell band $\widetilde E_m$
satisfies
\[
 \widetilde E_m'(k)=(E_j^{(n)})'(k+2\pi\nu/P)
\]
for a fixed pair $(j,\nu)$, with $j\ge1$ and $0\le\nu<\ell$,
on each nonsymmetric half-zone.
Figure~\ref{fig:folding} illustrates both constancy within a half-zone and
a change of label at its boundary.

Put
\begin{equation}\label{def:universalc}
 c\eqdef\frac{\sqrt{2R}}2.
\end{equation}
The small-period discriminant estimate controls these derivatives with
exponential cost $e^{2cp}$. Large-cell monotonicity and non-overlap then
separate all ordered eigenvalues away from the symmetry set. A bounded
perturbation changes each ordered eigenvalue by at most its norm.
Consequently, for the next approximant, every pair at distance at most
$\varepsilon$ lies in the same neighborhood of the symmetry momenta,
and its measure is bounded by a constant times $e^{cp_n}\sqrt\varepsilon$
when $\varepsilon>\|V_{n+1}-V_n\|_\infty$.

The continuum implementation requires two additional estimates. Each
fiber has infinitely many bands, so the near-collision neighborhood must
be independent of the band index. The momentum $D=-i\partial_x$ is
unbounded, so the averaging and perturbation arguments require form and
graph-norm control. The finite-dimensional folding and min--max comparison
already occur in \cite{HGS}; the argument here replaces matrix-size bounds
by a period-uniform form estimate and a sum over ordered eigenvalue
differences. It also proves a linear potential-error estimate on
$H^2(\R)\cap\Dom(X)$ by a common-domain graph-norm identity.

For the averaged momentum $\mathcal A_H(t)$ defined in
\eqref{eq:Av}, the resulting estimate is
\begin{equation}\label{eq:intromaster}
 \left\|\mathcal A_{H_{n+1}}(t)\psi-\tfrac12Q_{n+1}\psi\right\|^2
 \le C(R)C_{\psi,s}
 \left[p_{n+1}e^{cp_n}\sqrt\varepsilon+
                   \frac1{t^2\varepsilon^2}\right].
\end{equation}
It applies for $\psi\in\mathcal D_s$, $s>2$, and
$\varepsilon>\|V_{n+1}-V_n\|_\infty$ after the harmless normalization
of Section~\ref{sec:fold}; see Corollary~\ref{cor:EII}.
The remaining factor $p_{n+1}$ comes from the normalized measure
$(p_{n+1}/2\pi)\,\dd k$. The rate estimate and the linear dynamical
comparison determine common comparison times for successive approximants.
Their errors are summable and smaller than a finite-scale lower bound on
the velocity, yielding both convergence and nonvanishing.

\subsection{Organization}
Section~\ref{sec:periodic} fixes the Floquet conventions and discriminant
bound. Section~\ref{sec:fold} proves exact folding and the fixed-label
statement by ordering. Section~\ref{sec:badset} locates the near-collisions;
folded band widths are recorded there as a separate geometric consequence.
Section~\ref{sec:rate} proves the averaging estimate, and
Section~\ref{sec:main} proves Theorem~\ref{thm:main} and
Corollary~\ref{cor:moments}. Section~\ref{sec:remarks} discusses the state
class and open questions. The appendices supply the discriminant,
dynamical, and fiber-mass estimates used in the proofs.

\section{Preliminaries on periodic operators}\label{sec:periodic}
All inner products are conjugate-linear in the first argument. Unless
indicated otherwise, vector norms are those of the relevant Hilbert space
and operator norms are the corresponding induced norms. Potential norms
are written with the subscript $\infty$. For a real bounded $p$-periodic
potential, write
\[
 \T=\R/p\Z,\qquad \T^*=\R/(2\pi/p)\Z,\qquad
 \mathcal B=[-\pi/p,\pi/p].
\]
We use Lebesgue measure $\dd q$ on $\T$ and normalized measure
$\dd k/|\T^*|=(p/2\pi)\dd k$ on $\mathcal B$.

\begin{definition}[Floquet transform]
The unitary Floquet transform is
\[
 \begin{gathered}
 U\colon L^2(\R)\longrightarrow
 \int_{\T^*}^{\oplus}L^2(\T,\dd q)\,\frac{\dd k}{|\T^*|},\\
 (U\psi)(k,q)=\sum_{m\in\Z}e^{-ik(q+pm)}\psi(q+pm).
 \end{gathered}
\]
For a general $L^2$ input, the series converges in the direct-integral $L^2$
norm. Its inverse is
\[
 (U^*f)(q+pm)=\int_{\T^*}e^{ik(q+pm)}f(k,q)\,
 \frac{\dd k}{|\T^*|},\qquad 0\le q<p,
\]
with the corresponding $L^2$ interpretation.
\end{definition}

All integrals over $\T^*$ are understood over the chosen fundamental
interval $\mathcal B$; their fiber data obey the gauge identification below.
For $G=2\pi/p$, the extensions to the real $k$-axis satisfy the reciprocal-lattice gauge identities
\begin{equation}\label{eq:gauge}
 (U\psi)(k+G,q)=e^{-iGq}(U\psi)(k,q),\qquad
 H(k+G)=e^{-iGq}H(k)e^{iGq}.
\end{equation}
Thus fiber vectors and operators are gauge-periodic, whereas ordered
eigenvalues are $G$-periodic.

With $D=-i\partial_x$, one has
\begin{equation}\label{eq:intertwine}
 UDU^*=\int_{\T^*}^{\oplus}(D+k)\,\frac{\dd k}{|\T^*|},
\end{equation}
and
\begin{equation}\label{eq:direct}
 UHU^*=\int_{\T^*}^{\oplus}H(k)\,\frac{\dd k}{|\T^*|},\qquad
 H(k)=(D+k)^2+V.
\end{equation}
The identities are first obtained on $C_c^\infty(\R)$, where the differentiated
Floquet sums are finite on each cell, and then extended in the corresponding
graph and form norms. The common operator domain is $H^2(\T)$ and the common form domain is
$H^1(\T)$. Each real fiber is self-adjoint, bounded below by
$-\|V\|_\infty$, and has compact resolvent. The free fiber is
$H_0(k)=(D+k)^2$.

The family $z\mapsto(D+z)^2+V$ is analytic on the common domain
$H^2(\T)$: the lower-order perturbation $2zD+V+z^2$ has relative bound
zero with respect to $D^2$, and the action on each domain vector is
polynomial in $z$. Local analytic perturbation theory \cite[Chapter~VII]{Kat}
therefore applies to simple isolated eigenvalues and their projections.

For a displayed period $p$, let $\mathcal M_p(\lambda)$ carry
$(u(0),u'(0))$ to $(u(p),u'(p))$ for $-u''+Vu=\lambda u$.
The discriminant is $\Delta(\lambda)=\Delta_p(\lambda)
\eqdef\tr\mathcal M_p(\lambda)$, for the displayed, possibly nonminimal,
period. The monodromy determinant is one by constancy of the Wronskian.

\begin{proposition}[Floquet band properties]\label{prop:floquet}
Let $E_m(k)$ be the ordered eigenvalues of $H(k)$, counted with multiplicity.
For $k\notin(\pi/p)\Z$, let $P_m(k)$ be their rank-one spectral projections.
At the symmetry points, choose any orthogonal rank-one splitting of each
multiple eigenspace. These values on a null set do not affect direct integrals.
The following properties hold for every displayed period $p$, whether or
not it is minimal:
\begin{enumerate}[label=\textup{(F\arabic*)}]
\item\label{F1} $E_m(k)=E_m(-k)$.
\item\label{F2} Away from $(\pi/p)\Z$, the eigenvalues are simple. On
$(0,\pi/p)$, the function $E_m$ is real-analytic and
$(-1)^{m+1}E_m'(k)>0$. Its minimum on $\mathcal B$ is at $0$ for $m$ odd
and at $\pm\pi/p$ for $m$ even; its maximum is at the other symmetry point.
\item\label{F3} The bands do not overlap:
$\max_{\mathcal B}E_m\le\min_{\mathcal B}E_{m+1}$.
\item\label{F4} The Hill discriminant obeys
$\Delta(E_m(k))=2\cos(pk)$ for every real $k$.
\item\label{F5} For almost every $k$,
\begin{equation}\label{eq:GV}
 P_m(k)(D+k)P_m(k)=\tfrac12E_m'(k)P_m(k).
\end{equation}
\end{enumerate}
\end{proposition}
The first four facts are standard one-dimensional Floquet theory; see
\cite[Section 2]{Young} and \cite[Section XIII.16]{RSIV}. Identity~\eqref{eq:GV} is the
Hellmann--Feynman formula used by Asch and Knauf \cite{AsKn}.
On $0<k<\pi/p$, differentiating (F4) gives
$\Delta'(E_m(k))E_m'(k)=-2p\sin(pk)\ne0$; this supplies the strictness in (F2).

\begin{lemma}[Uniform discriminant bound]\label{lem:discriminant}
Put $C_2=2e^3$ and $c=\sqrt{2R}/2$. For every $R>0$, every real
bounded $p$-periodic potential with $\|V\|_\infty\le R$, and every
$\lambda\ge-R$,
\begin{equation}\label{eq:disc}
 |\Delta'(\lambda)|\le C_2p^2e^{2cp}.
\end{equation}
\end{lemma}
A related estimate on the spectrum is recorded in
\cite[Lemma 2.1]{Young}; see also \cite[Lemma 2.1]{Fil2015}.
Appendix~\ref{app:disc} proves the stated bound for all $\lambda\ge-R$,
with the prefactor and exponential rate specified separately.

\section{The geometry of period enlargement}\label{sec:fold}
The comparison will use the old operator in two representations. Its
small cell controls derivatives, while its large cell supplies ordered bands
that can be compared directly with the next approximant. The construction below fixes the small-cell label of each ordered
large-cell band on a nonsymmetric half-zone.

\subsection{Nested approximants and normalization}
Under Definition~\ref{def:class}, the continuous approximants are uniformly
Cauchy: the essential supremum of their difference equals its ordinary
supremum. Their uniform limit is a continuous representative of $V$.
We use that representative, which defines the same operator.

\begin{proposition}[{cf.\ \cite[Proposition 2.1]{HGS}}]\label{prop:telescoping}
Let $H\in\EC(\eta,\al)$ with approximants $V_n$. Set $Z_1=V_1$ and
$Z_k=V_k-V_{k-1}$ for $k\ge2$. Then each $Z_k$ is continuous and
$p_k$-periodic, $\sum_{k=1}^N Z_k$ converges uniformly to the chosen
representative of $V$, and there is $N_1$ such that
\[
 \|Z_k\|_\infty\le2p_k^{-\al}e^{-\eta p_{k-1}},\qquad k\ge N_1.
\]
\end{proposition}
\begin{proof}
The assertion for $Z_1$ is immediate. For $k\ge2$, both $V_k$ and $V_{k-1}$
are $p_k$-periodic, since $p_k/p_{k-1}\in\N$. The partial sum is $V_N$, which
converges uniformly to $V$. Choose $N_*$ so that
$\|V-V_n\|_\infty\le p_{n+1}^{-\al}e^{-\eta p_n}$ for $n\ge N_*$. Then,
for $k\ge N_*+1$,
\[
 \|Z_k\|_\infty\le p_{k+1}^{-\al}e^{-\eta p_k}
      +p_k^{-\al}e^{-\eta p_{k-1}}
 \le2p_k^{-\al}e^{-\eta p_{k-1}}.
\]
Take $N_1=N_*+1$.
\end{proof}

The following clipping argument permits us to assume, without changing $R$,
that $\|V_n\|_\infty\le R$ for every $n$. For $R>0$, define the \emph{clipping map}
\begin{equation}\label{eq:clipmap}
    \mathfrak c_R(t)\eqdef\max\{-R,\min\{t,R\}\},\qquad t\in\R.
\end{equation}
Thus $\mathfrak c_R$ is the metric projection of $\R$ onto the closed interval
$[-R,R]$. In particular, $\mathfrak c_R(t)=t$ for $|t|\le R$ and
\begin{equation}\label{eq:cliplip}
    |\mathfrak c_R(t)-\mathfrak c_R(s)|\le |t-s|,\qquad s,t\in\R,
\end{equation}
so the clipping map is $1$-Lipschitz. Replace each approximant by
\begin{equation}\label{eq:clip}
    \widehat V_n(x)\eqdef \mathfrak c_R(V_n(x)).
\end{equation}
Since $\|V\|_\infty\le R$, one has $\mathfrak c_R(V(x))=V(x)$ for a.e. $x$.
Consequently $\widehat V_n$ is continuous and $p_n$-periodic and, by
\eqref{eq:cliplip},
\begin{equation}
    \|V-\widehat V_n\|_\infty\leq \|V-V_n\|_\infty.
\end{equation}
Relabeling $\widehat V_n$ as $V_n$ preserves the hypothesis of Definition~\ref{def:class}. For $j\ge2$ set
\begin{equation}\label{eq:gammadef}
    \gamma_j\eqdef\|V_j-V_{j-1}\|_\infty.
\end{equation}
The clipping map does not increase successive differences, so the
argument of Proposition~\ref{prop:telescoping} gives, after discarding finitely many terms,
\begin{equation}\label{eq:gamma}
    \gamma_{n+1}\leq 2p_{n+1}^{-\al}e^{-\eta p_n}.
\end{equation}
After discarding finitely many further terms, we also assume $p_n\ge\pi$.
Reindexing this tail preserves the period-ratio condition and the limit in
Definition~\ref{def:class}.

For each $j\in\N$, set
\begin{equation}\label{eq:Hj}
    H_j\eqdef-\partial_x^2+V_j,
    \qquad \ell_j\eqdef\frac{p_{j+1}}{p_j}\in\N,\qquad \ell_j\ge2.
\end{equation}
Fix $n\in\N$ and write $\ell\eqdef\ell_n$. We make the following convention about the tilde notation:
\begin{equation}\label{eq:Htilde}
    \widetilde H_n\eqdef -\partial_x^2+V_n
    \qquad \text{regarded as a $p_{n+1}$-periodic operator.}
\end{equation}

We also write 
\begin{equation}
    \T_n \eqdef \R/p_n\Z, \qquad \T^\ast_{n}\eqdef\R\left/\frac{2\pi}{p_n}\Z\right., \qquad \mathcal{B}_n \eqdef \left[-\frac{\pi}{p_n},\frac{\pi}{p_n}\right].
\end{equation}
For the $p_j$-periodic operator $H_j$, we write $U_{p_j}$ for the Floquet transform
of Section 2 and denote its ordered fiber eigenvalues and spectral projections by
$E_m^{(j)}(k)$ and $P_m^{(j)}(k)$, respectively.

This means that $\widetilde H_n$ has the fiber representation $\widetilde H_n(k)=(D+k)^2+V_n$ on $H^2(\T_{n+1})$, for real $k$, with ordered eigenvalues $\{\widetilde E_m(k)\}_{m=1}^\infty$.
These eigenvalues are even and $2\pi/p_{n+1}$-periodic; $\mathcal B_{n+1}$
is their chosen fundamental interval.

For the exact folding discussion fix $n$ and abbreviate
\[
 p=p_n,\quad {P}=p_{n+1}=\ell p,\quad
 \widetilde H(k)=\widetilde H_n(k),\quad
 H_p(\kappa)=H_n(\kappa).
\]
Write $D_L=-i\partial_x$ on the circle of length $L$. Its second-order
periodic domain is
\[
 {H^2_{\rm per}(0,L)}=\{u\in H^2(0,L):u(L)=u(0),\ u'(L)=u'(0)\};
\]
the form domain has the single periodic trace condition in $H^1(0,L)$.
The operators $(D_L+k)^2+V_n$ have these common domains. Boundedness of $V_n$
and the Fourier basis give compact resolvent and a discrete ordered
spectrum. Throughout this section, ${P}$ denotes the large period and
$\mathsf S_\nu$ a translation-sector projection. Whenever
$\widetilde E_m(k)$ is simple, $\Pi_m(k)$ denotes its rank-one spectral
projection for $\widetilde H(k)$. 

\subsection{Translation sectors and exact folding}
The next two lemmas give a continuum, domain-compatible version of the
period-enlargement decomposition used in the proof of
\cite[Theorem 3.2]{HGS}.

\begin{lemma}[Translation sectors and their domains]\label{fold:sectors}
Let $\omega=e^{2\pi i/\ell}$ and $(\tau u)(x)=u(x+p)$, with the argument
taken modulo ${P}$. The operators
\begin{equation}\label{fold:eq-sectors}
 \mathsf S_\nu=\frac1\ell\sum_{a=0}^{\ell-1}\omega^{-\nu a}\tau^a,
 \qquad 0\le \nu<\ell,
\end{equation}
are orthogonal projections satisfying
\[
 \sum_{\nu=0}^{\ell-1}\mathsf S_\nu=I,\qquad
 \tau\mathsf S_\nu=\omega^\nu\mathsf S_\nu,\qquad
 \Ran\mathsf S_\nu=K_\nu:=\{u:\tau u=\omega^\nu u\}.
\]
They preserve the operator and form domains of $\widetilde H(k)$ and
commute with its resolvent, also for complex $k$ wherever the resolvent
exists.
\end{lemma}
\begin{proof}
Translation preserves measure, so $\tau^*=\tau^{-1}$ and $\tau^\ell=I$.
Reindexing $a+1$ modulo $\ell$ gives $\tau\mathsf S_\nu=\omega^\nu\mathsf S_\nu$;
reindexing $-a$ gives $\mathsf S_\nu^*=\mathsf S_\nu$. In the product of the
two finite sums, the identity
\[
 \sum_{a=0}^{\ell-1}\omega^{(s-\nu)a}
 =\begin{cases}\ell,&\nu=s,\\0,&\nu\ne s\end{cases}
\]
gives $\mathsf S_\nu\mathsf S_s=\delta_{\nu s}\mathsf S_\nu$. Summing first over
$\nu$ gives $\sum_\nu\mathsf S_\nu=I$, proving the complete orthogonal
sector decomposition.

Translation preserves periodic Sobolev functions and commutes with their
weak derivatives. It also commutes with $V_n$, since $V_n(x+p)=V_n(x)$ almost
everywhere. Thus $\tau{H^2_{\rm per}(0,{P})}={H^2_{\rm per}(0,{P})}$ and
$\tau\widetilde H(k)u=\widetilde H(k)\tau u$ for $u\in{H^2_{\rm per}(0,{P})}$.
If $u=(z-\widetilde H(k))^{-1}f$, uniqueness of the resolvent equation
applied to
\[
 (z-\widetilde H(k))\tau u=\tau f
\]
shows that the resolvent commutes with $\tau$. Each $\mathsf S_\nu$ is a finite polynomial
in $\tau$, so all the asserted commutation and domain properties follow.
\end{proof}

\begin{lemma}[Folding identity]\label{lem:fold}
For $\vartheta_{\nu}=2\pi \nu/{P}$, the map
\[
 (\iota_\nu f)(x)=\ell^{-1/2}e^{i\vartheta_{\nu}x}f(x\bmod p)
\]
is an isometry from $L^2(\T_n)$ onto $K_\nu$, and
\begin{equation}\label{fold:eq-intertwining}
 \widetilde H(k)\iota_\nu=\iota_\nu H_p(k+\vartheta_{\nu}),\qquad
 \iota_\nu^*\iota_\nu=I,\quad \iota_\nu\iota_\nu^*=\mathsf S_\nu.
\end{equation}
The first identity respects both operator and form domains. In particular,
counting multiplicities,
\begin{equation}\label{eq:fold-spectrum}
 \sigma(\widetilde H_n(k))
 =\bigcup_{\nu=0}^{\ell-1}
   \{E_j^{(n)}(k+2\pi \nu/{P}):j\ge1\}.
\end{equation}
\end{lemma}
\begin{proof}
The restriction of the inverse map to the first small cell is
\[
 W_\nu:K_\nu\longrightarrow L^2(\T_n),\qquad
 (W_{\nu}u)(x)=\ell^{1/2}e^{-i\vartheta_{\nu}x}u(x),\quad 0\le x<p.
\]
The phase cancels the character $\omega^\nu$ under translation by $p$.
Moreover $|u|$ is $p$-periodic, and therefore
\[
 \|W_{\nu}u\|_{L^2(0,p)}^2
 =\ell\int_0^p|u|^2=\int_0^{P}|u|^2.
\]
The phase-modulated periodic extension $\iota_\nu$ is its inverse, proving
unitarity. The same phase relation holds for the first derivative at a cell
boundary. Multiplication by $e^{-i\vartheta_{\nu}x}$ turns both sector trace
conditions into periodic trace conditions. Hence $W_\nu$ carries
${H^2_{\rm per}(0,{P})}\cap K_\nu$ onto ${H^2_{\rm per}(0,p)}$ and has the analogous property
on the form domain.

On these domains,
\[
 D_{P}(e^{i\vartheta_{\nu}x}f)=e^{i\vartheta_{\nu}x}(D_p+\vartheta_{\nu})f,
 \qquad
 W_\nu\widetilde H(k)W_\nu^{-1}=(D_p+k+\vartheta_{\nu})^2+V_n.
\]
Lemma~\ref{fold:sectors} now gives the orthogonal direct sum over $\nu$,
including multiplicity, and proves \eqref{eq:fold-spectrum}.
\end{proof}
The sector decomposition identifies the translated small-cell fibers. Write
\[
 \kappa_\nu(k)=k+2\pi \nu/{P},\qquad
 \varphi_{j,\nu}(k)=E_j^{(n)}(\kappa_\nu(k)).
\]
The eigenvalues $\widetilde E_m(k)$ are obtained by ordering these numbers.

\subsection{Where an ordered branch can change its label}
Proposition~\ref{prop:floquet} applies to the old operator with displayed period $P$, even when this period is not minimal.

\begin{lemma}[Simplicity and uniqueness of a folded pair]\label{fold:simplicity}
For real $k\notin(\pi/{P})\Z$, every eigenvalue of $\widetilde H(k)$ is
simple. Its occurrence in \eqref{eq:fold-spectrum} is represented by exactly
one pair $(j,\nu)$.
\end{lemma}
\begin{proof}
An eigenfunction $v\in{H^2_{\rm per}(0,{P})}$ of energy $E$ gives
$u(x)=e^{ikx}v(x)$ satisfying
\[
 -u''+V_n u=Eu,\qquad
 \binom{u({P})}{u'({P})}=e^{ik P}\binom{u(0)}{u'(0)}.
\]
For bounded measurable $V_n$, the corresponding first-order Volterra system
has a unique solution for each initial pair. Its monodromy is real for
real $E$ and has determinant one, since the Wronskian is constant.
Its eigenvalues are therefore $e^{ik P}$ and $e^{-ik P}$. They are distinct
off $(\pi/{P})\Z$, so the eigenspace corresponding to $e^{ik P}$ is one-dimensional. The map from
a fiber eigenfunction to its initial pair is injective; self-adjointness
then gives spectral multiplicity one. Minimality of ${P}$ is not used.

If two folded pairs had the same energy, \eqref{fold:eq-intertwining}
would give two independent large-cell eigenvectors. Different sectors are
orthogonal; equal sectors preserve the small-cell multiplicity under
$\iota_\nu$. Both alternatives contradict simplicity.
\end{proof}
The shifted small momentum also avoids its own symmetry set:
\begin{equation}\label{fold:eq-avoid-small}
 k+\frac{2\pi \nu}{{P}}=\frac{\pi a}{p}\quad(a\in\Z)
 \quad\Longrightarrow\quad
 k=\frac\pi {P}(\ell a-2\nu)\in\frac\pi {P}\Z.
\end{equation}

The following form comparison uses the min--max principle
\cite[Section XIII.1]{RSIV}.
\begin{lemma}[Continuity of the ordered bands]\label{lem:lip}
On a compact real interval $I$, each $\widetilde E_m$ is Lipschitz. If
$\Lambda_m=\max_I\widetilde E_m$, then
\begin{equation}\label{eq:lip}
 |\widetilde E_m(k')-\widetilde E_m(k)|
 \le2|k'-k|\sqrt{\Lambda_m+R}+|k'-k|^2,\qquad k,k'\in I.
\end{equation}
The same assertion holds for each ordered small-cell band.
\end{lemma}
\begin{proof}
Finiteness of a common upper bound does not require continuity. Put
$K_I=\max_{k\in I}|k|$ and use the $m$-dimensional span of the Fourier
modes ${P}^{-1/2}e^{2\pi iax/{P}}$, $a=0,\ldots,m-1$, as a trial space.
Min--max gives
\[
 \widetilde E_m(k)\le
 \left(K_I+\frac{2\pi(m-1)}{P}\right)^2+R\eqdef L_{m,I}.
\]
For a unit vector $u$ in the span of the first $m$ eigenvectors at $k$,
\begin{equation}\label{eq:momu}
 \|(D_{P}+k)u\|^2\le\widetilde E_m(k)+R.
\end{equation}
Expanding the common quadratic form at $k'$ yields
\[
 \widetilde{\mathfrak h}_{k'}[u]
 =\widetilde{\mathfrak h}_k[u]
 +2(k'-k)\operatorname{Re}\langle u,(D_{P}+k)u\rangle
 +(k'-k)^2.
\]
Cauchy--Schwarz and min--max imply
\[
 \widetilde E_m(k')\le\widetilde E_m(k)
 +2|k'-k|\sqrt{\widetilde E_m(k)+R}+|k'-k|^2.
\]
Exchange $k$ and $k'$ and first use the trial-space bound $L_{m,I}$.
This proves Lipschitz continuity, since
$|k'-k|^2\le\operatorname{diam}(I)|k'-k|$. The continuous function now has
an attained maximum $\Lambda_m$, giving \eqref{eq:lip} as stated.
The identical argument with cell length $p$ proves the small-cell assertion.
In particular, the ordered bands are absolutely continuous up to symmetry
endpoints and their almost-everywhere derivatives may be integrated.
\end{proof}

\medskip\noindent The fixed-label assertion now follows by ordering.
\begin{proposition}[Fixed small-cell labels]\label{prop:fixedlabels}
On each component $J$ of $\R\setminus(\pi/{P})\Z$ and for every $m\ge1$,
a unique pair $(j,\nu)$, independent of $k\in J$, satisfies
\begin{equation}\label{eq:fixedpair}
 \widetilde E_m(k)=E_j^{(n)}(k+\vartheta_{\nu}).
\end{equation}
This eigenvalue is real-analytic on $J$, and
\begin{equation}\label{eq:foldderivative}
 \widetilde E_m'(k)=(E_j^{(n)})'(k+\vartheta_{\nu}).
\end{equation}
\end{proposition}
\begin{proof}
Every translated small-cell branch $\varphi_{j,\nu}$ is continuous by
Lemma~\ref{lem:lip}. Lemma~\ref{fold:simplicity} says that two different
pairs never give the same energy on $J$. Their difference is therefore a
continuous nonzero function on an interval and has a fixed sign.

Fix $k_0\in J$ and the unique pair whose energy is
$\widetilde E_m(k_0)$. Exactly $m-1$ other pairs lie below it at $k_0$.
This number is finite because the large-cell operator has compact
resolvent and is bounded below. Pairwise order is unchanged throughout
$J$, so precisely those same pairs, and no others, lie below the chosen
pair at every $k\in J$. Its ordered index remains $m$.
The avoidance relation \eqref{fold:eq-avoid-small} places $k+\vartheta_{\nu}$
away from the small-cell symmetry set. Proposition~\ref{prop:floquet}
then gives real analyticity. Differentiation proves
\eqref{eq:foldderivative}, since $\dd(k+\vartheta_{\nu})/\dd k=1$.
\end{proof}
\begin{localremark}[The associated eigenprojection]
Intertwining in Lemma~\ref{lem:fold} and simplicity also give
\begin{equation}\label{eq:foldprojections}
 \Pi_m(k)=\iota_{\nu}P_j^{(n)}(k+\vartheta_{\nu})\iota_\nu^*\qquad(k\in J).
\end{equation}
Both sides are the orthogonal projection onto the same one-dimensional
eigenspace. Analytic perturbation makes this map real-analytic on $J$.
The derivative estimate below uses the eigenvalue identity
\eqref{eq:fixedpair}; the projection formula records the corresponding
spectral subspace.
\end{localremark}

\begin{example}\label{ex:switch}
Take $V_n=0$, $p_n=p\ge\pi$, and $\ell=3$, so ${P}=p_{n+1}=3p$.
On $\mathcal B_{n+1}=[-\pi/(3p),\pi/(3p)]$, the folded eigenvalues are
$\{(k+2\pi a/(3p))^2:a\in\Z\}$. Since
$E_1^{(n)}(\kappa)=\dist(\kappa,(2\pi/p)\Z)^2$,
\[
 \widetilde E_2(k)=\left(\frac{2\pi}{3p}-|k|\right)^2
 =\begin{cases}
   \varphi_{1,1}(k),&-\pi/(3p)<k<0,\\
   \varphi_{1,2}(k),&0<k<\pi/(3p).
  \end{cases}
\]
Thus the pair $(j,\nu)$ changes at the large-period symmetry point $k=0$,
where the one-sided derivatives are $4\pi/(3p)$ and $-4\pi/(3p)$.
This is a crossing of different translation sectors. A corner need not
always change the pair: for $\ell=2$, the second folded band is
$E_1^{(n)}(k+\pi/p)=(\pi/p-|k|)^2$ throughout the large zone, and its
corner comes from the small-cell symmetry point with the same label $(1,1)$.
\end{example}

\subsection{The smaller-period derivative on the folded band}
The following band-velocity and quadratic-separation estimates are due to
Young. We follow his proof, using the discriminant constants of
Lemma~\ref{lem:discriminant}, before applying the folding identity.
\begin{proposition}[{\cite[Lemma 2.2]{Young}}]\label{prop:vel}
For every $p$-periodic $V$ with $\|V\|_\infty\leq R$, every $m\geq1$, and every $k\in\R\setminus(\pi/p)\Z$ (using the even, $2\pi/p$-periodic extension of $E_m$),
\begin{equation}\label{eq:26}
    \left|\frac{\dd E_m}{\dd k}(k)\right|\geq \frac{2|\sin(pk)|}{C_2p}e^{-2cp}.
\end{equation}
In particular, the estimate holds for a.e. $k\in\mathcal B$. Moreover, for $0\le k_1\le k_2\le\pi/p$,
\begin{equation}\label{eq:27}
    |E_m(k_2)-E_m(k_1)|\geq \frac{e^{-2cp}}{\pi C_2}(k_2-k_1)^2.
\end{equation}
\end{proposition}
\begin{proof}
The differentiated discriminant identity gives \eqref{eq:26}. For
$0\le k_1\le k_2\le\pi/p$, monotonicity and absolute continuity give
\[
 |E_m(k_2)-E_m(k_1)|\ge\frac{2e^{-2cp}}{C_2p}
       \int_{k_1}^{k_2}\sin(ps)\,\dd s.
\]
Among intervals of length $h$ in $[0,\pi/p]$, this integral is minimized
when the interval abuts an endpoint; its minimum is $(1-\cos(ph))/p$. Since
$1-\cos(ph)\ge2p^2h^2/\pi^2$, the resulting bound
$4e^{-2cp}h^2/(C_2\pi^2)$ implies \eqref{eq:27}.
\end{proof}

For $k\in\mathcal B_{n+1}$, put ${P}=p_{n+1}$ and define the distance
to the large-period symmetry set by
\begin{equation}\label{eq:dk}
    d_{P}(k)\eqdef\dist\!\left(k,\frac{\pi}{p_{n+1}}\Z\right)
    =\min\left\{|k|,\frac{\pi}{p_{n+1}}-|k|\right\}.
\end{equation}

The next estimate is the physical-quasimomentum counterpart of the
folded-slope calculation in \cite[proof of Theorem 3.2]{HGS}.
\begin{lemma}\label{lem:folded}
For every $m\ge1$ and a.e.\ $k\in\Bz_{n+1}$,
\begin{equation}\label{eq:folded}
\Big|\frac{\dd\widetilde E_m}{\dd k}(k)\Big|\ \ge\ \frac{4}{\pi C_2}\,e^{-2cp_n}\,d_{P}(k).
\end{equation}
\end{lemma}

\begin{proof}
By Proposition~\ref{prop:fixedlabels}, off the symmetry set $\widetilde E_m$ coincides with a
single branch $\varphi_{j,\nu}$; fix such a $k$ and set $\kappa\eqdef\kappa_\nu(k)$.
Since $\kappa_\nu'(k)=1$,
$\tfrac{\dd\widetilde E_m}{\dd k}(k)=\tfrac{\dd E^{(n)}_j}{\dd\kappa}(\kappa)$, with no
chain-rule factor, and Proposition~\ref{prop:vel} at period $p_n$ gives
\begin{equation}\label{eq:step1}
\Big|\frac{\dd\widetilde E_m}{\dd k}(k)\Big|\ \ge\ \frac{2}{C_2p_n}e^{-2cp_n}\big|\sin(p_n\kappa)\big| .
\end{equation}
Put $x\eqdef p_{n+1}k\in[-\pi,\pi]$, so that $p_n\kappa=(x+2\pi\nu)/\ell$. Using successively
$\dist(y/\ell,\pi\Z)=\tfrac1\ell\dist(y,\ell\pi\Z)$,
$\ell\pi\Z\subseteq\pi\Z$, and $2\pi\nu\in\pi\Z$, we obtain
\[
\begin{aligned}
\dist\Big(\frac{x+2\pi\nu}{\ell},\pi\Z\Big)
&=\frac1\ell\dist(x+2\pi\nu,\ell\pi\Z)\\
&\ge\frac1\ell\dist(x+2\pi\nu,\pi\Z)
=\frac1\ell\dist(x,\pi\Z).
\end{aligned}
\]
Since $|\sin u|\ge\tfrac2\pi\dist(u,\pi\Z)$ and, for $k\in\Bz_{n+1}$,
$\dist(p_{n+1}k,\pi\Z)=\min\{|x|,\pi-|x|\}=p_{n+1}d_{P}(k)$, we get
$|\sin(p_n\kappa)|\ge\tfrac{2}{\pi\ell}p_{n+1}d_{P}(k)$. Substituting into
\eqref{eq:step1} and using $p_{n+1}=\ell p_n$ gives \eqref{eq:folded}.
\end{proof}

\begin{remark}[Quasimomentum coordinates]\label{rem:cancel}
The calculation above is in physical quasimomentum. If instead
\[
 \theta=\frac{p_{n+1}k}{2\pi},\qquad
 \theta_n=\frac{p_n\kappa_\nu(k)}{2\pi},
\]
then $\theta_n=(\theta+\nu)/\ell$ and the chain rule has a factor $\ell^{-1}$.
Thus translation versus dilation depends on the coordinate convention, not
on whether the model is continuous or discrete. Comparisons with the
unit-length quasimomentum zone in \cite{HGS} must also convert the measure.
The cancellation in \eqref{eq:folded} is valid in the physical-$k$ convention;
the improvement of the exponential scale from $p_{n+1}$ to $p_n$ is
independent of this convention.
\end{remark}

\section{Spectral consequences and near-collision geometry}\label{sec:badset}
The fixed-label identity transfers the old-period slope bound to the
ordered large-cell bands. Non-overlap converts that bound into separation
away from symmetry momenta. Min--max then locates the near-collisions of
the next approximant.

\subsection{Eigenvalue stability in the common large cell}
\begin{proposition}[{cf.\ \cite[Proposition 3.1]{HGS}}]\label{prop:minmax}
For every $m\ge1$ and $k\in\Bz_{n+1}$,
$\ \big|E^{(n+1)}_m(k)-\widetilde E_m(k)\big|\le\gamma_{n+1}$.
\end{proposition}

\begin{proof}
Both $H_{n+1}(k)$ and $\widetilde H_n(k)$ are self-adjoint on $H^2(\T_{n+1})$, bounded below
by $-R$, and are bounded perturbations of $(D+k)^2$ on a torus, hence have
compact resolvent and purely discrete spectrum accumulating only at $+\infty$.
The Courant--Fischer principle \cite[Section XIII.1]{RSIV} gives
\[
\begin{gathered}
E^{(n+1)}_m(k)=\min_{\substack{S\subset H^1(\T_{n+1})\\\dim S=m}}\ \max_{\substack{u\in S\\\|u\|=1}}\mathfrak h_k[u],\\
\mathfrak h_k[u]=\|(D+k)u\|^2+\int_{\T_{n+1}}V_{n+1}|u|^2,
\end{gathered}
\]
and likewise for $\widetilde E_m(k)$ with $V_n$ in place of $V_{n+1}$. The two forms have
the same form domain $H^1(\T_{n+1})$ and
$|\mathfrak h_k[u]-\widetilde{\mathfrak h}_k[u]|\le\gamma_{n+1}\|u\|^2$. Taking
$\max$ over unit vectors of a fixed $S$ and then $\min$ over $S$ preserves this
two-sided bound.
\end{proof}

This comparison does not require the new potential $V_{n+1}$ to commute
with translation by $p_n$. The translation-sector statements apply only
to the old potential read in the large cell; min--max is the bridge to the
new operator.

\subsection{Separation away from the symmetry momenta}
\begin{proposition}[Separation of folded eigenvalues]\label{prop:foldgeometry}
Put $p=p_n$, ${P}=p_{n+1}$ and $d_{P}(k)=\dist(k,(\pi/{P})\Z)$.
For every $m<m'$ and $k\in\mathcal B_{n+1}$,
\begin{equation}\label{eq:foldseparation}
 \widetilde E_{m'}(k)-\widetilde E_m(k)
 \ge\frac{2}{\pi C_2}e^{-2cp}d_{P}(k)^2.
\end{equation}
Consequently,
\begin{equation}\label{eq:separation-stable}
 E_{m'}^{(n+1)}(k)-E_m^{(n+1)}(k)
 \ge\frac{2}{\pi C_2}e^{-2cp}d_{P}(k)^2-2\gamma_{n+1}.
\end{equation}
\end{proposition}
\begin{proof}
Let $L_*=\pi/{P}$. On $[0,L_*]$, $d_{P}(s)=\min\{s,L_*-s\}$.
If $z\in[0,L_*]$ is the distance from $k$ to the minimizing endpoint
of a band along its monotone half-zone, then
\[
 \begin{aligned}
 \int_0^z\min\{s,L_*-s\}\,\dd s
 &=\begin{cases}
 z^2/2,&z\le L_*/2,\\
 L_*z-z^2/2-L_*^2/4,&z\ge L_*/2,
 \end{cases}\\
 &\ge\tfrac12\min\{z,L_*-z\}^2.
 \end{aligned}
\]
For the second case, writing $u=L_*-z\le L_*/2$ makes the middle
expression $L_*^2/4-u^2/2\ge u^2/2$.
Monotonicity, absolute continuity, and Lemma~\ref{lem:folded} give
\begin{equation}\label{eq:quad}
 \widetilde E_m(k)-\min\widetilde E_m
 \ge\frac{2}{\pi C_2}e^{-2cp}d_{P}(k)^2.
\end{equation}
Evenness treats the negative half-zone. For $m<m'$, non-overlap gives
$\widetilde E_m(k)\le\max\widetilde E_m\le\min\widetilde E_{m'}$;
apply \eqref{eq:quad} to the upper band. Subtracting the two min--max
errors from Proposition~\ref{prop:minmax} proves
\eqref{eq:separation-stable}.
\end{proof}

\begin{remark}[A separate band-width consequence]\label{rem:spectral-scope}
The same slope estimate also gives a useful width bound:
\begin{equation}\label{eq:foldwidth}
 \max\widetilde E_m-\min\widetilde E_m
 =\int_0^{\pi/{P}}|\widetilde E_m'(s)|\,\dd s
 \ge\frac{\pi}{C_2{P}^2}e^{-2cp},
\end{equation}
since $\int_0^{\pi/{P}}d_{P}(s)\,\dd s=\pi^2/(4{P}^2)$.
Thus
\begin{equation}\label{eq:width-stable}
 \max E_m^{(n+1)}-\min E_m^{(n+1)}
 \ge\frac{\pi}{C_2{P}^2}e^{-2cp}-2\gamma_{n+1}.
\end{equation}
The extrema are over $\mathcal B_{n+1}$. These width bounds are geometric
byproducts; the transport proof uses \eqref{eq:foldseparation}.
Exact folding leaves $\sigma(H_n)$ unchanged, and the additional band
edges may be touching points. Width control alone does not imply gap
opening, Cantor spectrum, or homogeneity.
\end{remark}

\subsection{All near-collisions lie in one symmetry neighborhood}
For any ordered band family, the existence of a pair at distance at most
$\e$ is equivalent to the existence of an adjacent such pair. Hence its bad set is
\[
 B_\e=\bigcup_{m\ge1}
 \{k\in\mathcal B:E_{m+1}(k)-E_m(k)\le\e\}.
\]
Each set in this union is closed by Lemma~\ref{lem:lip}. Thus $B_\e$ is
Borel, both for a single period and for the approximating fibers below.
The following theorem implements the two-scale bad-set argument of
\cite[Theorem 3.2]{HGS} for continuum fibers, using the uniform separation
estimate proved above.
\begin{theorem}\label{thm:badset}
Let $\|V_n\|_\infty,\|V_{n+1}\|_\infty\le R$ and
$p_{n+1}=\ell p_n$, $\ell\ge2$. For $\e>0$, put
\[
 B^{(n+1)}_\e=\{k\in\mathcal B_{n+1}:\exists m\ne m',\
       |E_m^{(n+1)}(k)-E_{m'}^{(n+1)}(k)|\le\e\},
\]
and
\[
 r_\e=\left(\frac{\pi C_2}{2}\right)^{1/2}
 e^{cp_n}\sqrt{\e+2\gamma_{n+1}}.
\]
Then
\begin{equation}\label{eq:collision-location}
 B^{(n+1)}_\e\subseteq
 \{k\in\mathcal B_{n+1}:d_{P}(k)\le r_\e\},\qquad {P}=p_{n+1},
\end{equation}
and
\begin{equation}\label{eq:collision-general}
 \mes(B^{(n+1)}_\e)\le\min\{2\pi/{P},4r_\e\}.
\end{equation}
In particular, if $\gamma_{n+1}<\e$, then
\[
 \mes(B^{(n+1)}_\e)
 \le5\sqrt{C_2\pi}\,e^{cp_n}\sqrt\e.
\]
\end{theorem}
\begin{proof}
A colliding pair, reordered as $m<m'$, satisfies
\[
 \frac{2}{\pi C_2}e^{-2cp_n}d_{P}(k)^2
 \le E_{m'}^{(n+1)}(k)-E_m^{(n+1)}(k)+2\gamma_{n+1}
 \le\e+2\gamma_{n+1}
\]
by \eqref{eq:separation-stable}. This proves
\eqref{eq:collision-location}. If $r_\e\le\pi/(2{P})$, the indicated
set consists of two intervals adjacent to zero and one adjacent to each
endpoint $\pm\pi/{P}$, of total length $4r_\e$. If
$r_\e>\pi/(2{P})$, the full zone has length $2\pi/{P}<4r_\e$.
This proves \eqref{eq:collision-general} in both cases. Finally,
$\gamma_{n+1}<\e$ implies
\[
 4r_\e\le4\sqrt{3/2}\sqrt{\pi C_2}\,e^{cp_n}\sqrt\e
 \le5\sqrt{\pi C_2}\,e^{cp_n}\sqrt\e.
\]
The neighborhood is independent of the indices, so no infinite sum of
separate measure bounds is needed.
\end{proof}
Here $\mes$ is unnormalized Lebesgue measure. For the probability measure
in the direct integral, the last bound, under $\gamma_{n+1}<\e$, becomes
\begin{equation}\label{eq:normalizedbad}
 \int_{B^{(n+1)}_\e}\frac{\dd k}{|\T_{n+1}^*|}
 \le\frac{5\sqrt{C_2\pi}}{2\pi}
       p_{n+1}e^{cp_n}\sqrt\e.
\end{equation}
The larger period enters through the measure; the smaller period entered
through the derivative. Both factors have a specific geometric origin.

For comparison, the one-period estimate used below is the following.
\begin{proposition}[{cf.\ \cite[Lemma 2.3]{Young}}]\label{prop:badY}
If
\begin{equation}
    B_\e\eqdef\left\{k\in\mathcal B\colon \exists\,m\neq m',\ |E_m(k)-E_{m'}(k)|\leq\e\right\},
\end{equation}
then
\begin{equation}
    \mes(B_\e)\leq 4\sqrt{C_2\pi}\,e^{cp}\sqrt\e.
\end{equation}
\end{proposition}

Its proof with the constants of Lemma~\ref{lem:discriminant} is included
in Appendix~\ref{app:badone}.
\section{A periodic rate estimate with a free parameter}\label{sec:rate}
The periodic convergence argument of \cite[Theorem 3.1]{Young} separates
quasimomenta with small spectral gaps from their complement. We retain the
separation parameter explicitly, so that the same time can be used at two
consecutive periods. On the bad set a quadratic-form estimate suffices; on the
good set, ordering the eigenvalues allows summation in the index difference.
Neither estimate requires the full resolvent-weight sum.

\subsection{The periodic velocity and its form bounds}
For a $p$-periodic operator $H$, define
\begin{equation}\label{eq:Av}
    \mathcal A_H(t)\psi\eqdef\frac{1}{t}\int_0^tD(r)\psi\,\dd r,\qquad D(r)\eqdef e^{irH}De^{-irH}.
\end{equation}
The periodic asymptotic velocity is the following standard result.
Asch and Knauf formulate the hypotheses of \cite[Theorem 2.3]{AsKn}
for $d\ge2$ and use the kinetic energy $-\tfrac12\Delta$ when $\hbar=1$.
We use the bounded one-dimensional version below, with kinetic energy
$-\Delta$; the proof records the required form-domain argument.

\begin{theorem}[{Asch--Knauf; cf.\ \cite[Theorem 2.3]{AsKn} and \cite[Section 3]{Young}}]\label{thm:AK}
Let $H$ be $p$-periodic with bounded potential. For $\psi\in H^1(\R)\cap\Dom(X)$,
\begin{equation}
    \lim_{t\longrightarrow\infty}\frac1tX_H(t)\psi=Q_H\psi,
\end{equation}
where the limit is in $L^2(\R)$ and
\begin{equation}
    Q_H\psi\eqdef U^*\left(\int_{\T^\ast}^{\oplus}\sum_{m\geq1}\frac{\dd E_m(k)}{\dd k}P_m(k)\frac{\dd k}{|\T^\ast|}\right)U\psi.
\end{equation}
\end{theorem}

\begin{proof}
We record the form-domain justification in our normalization. Off the
symmetry points, put $A=H(k)+2R$ and $K=D+k$. The inequality $K^2\le A$
as forms gives
\[
 \sup_{t>0}\left\|\frac1t\int_0^t
 e^{i\tau H(k)}Ke^{-i\tau H(k)}v\,\dd\tau\right\|
 \le\|A^{1/2}v\|.
\]
For a vector in one eigenspace, expansion of $Kv$ in the complete
orthonormal eigenbasis reduces convergence to the scalar averages
$t^{-1}\int_0^t e^{i\tau(E_m-E_j)}\,\dd\tau$.
Their moduli are at most one, so dominated convergence in $\ell^2$
retains only $P_jKP_jv$. Finite spectral sums and then the tails
$\|A^{1/2}P_{>N}v\|\to0$ extend the result to $\Dom(A^{1/2})$.
The limiting diagonal sum has the same form bound. Dominated convergence
in the direct integral, Hellmann--Feynman, and the Heisenberg identity
in Lemma~\ref{lem:BHeis} finish the proof.
\end{proof}

For the approximants $H_n$ we use the abbreviation
\begin{equation}\label{eq:Qn}
    Q_n\eqdef Q_{H_n}.
\end{equation}

By \eqref{eq:GV}, $Q_H$ is twice the direct integral of
$\sum_m P_m(k)(D+k)P_m(k)$. Accordingly,
\begin{equation}\label{eq:33}
    X_H(t)\psi=X\psi+2\int_0^tD(r)\psi\,\dd r,
    \qquad \psi\in H^1(\R)\cap\Dom(X),
\end{equation}
as proved in Lemma~\ref{lem:BHeis}. The factor $2$ comes from
$i[H,X]=2D$ for our normalization $H=-\partial_x^2+V$, rather than
$-\tfrac12\partial_x^2+V$. The quantitative estimates below are therefore
written for $\mathcal A_H(t)\psi-\frac12Q_H\psi$.

\begin{proposition}[{cf.\ \cite[Lemma 3.3]{Young}}]\label{prop:M1}
If $\|V\|_\infty\leq R$, then
\begin{equation}
    \left\|(H(k)+2R)^{1/2}(D+k)(H(k)+2R)^{-1}\right\|^2\leq M_1
\end{equation}
uniformly in $k$ and in $p$; one may take $M_1=19$, independently of $R$.
We shall also use
\begin{equation}\label{eq:Dres}
    \left\|(D+k)(H(k)+2R)^{-1/2}\right\|\leq1.
\end{equation}
\end{proposition}
\begin{proof}
Put $A=H(k)+2R$, $K_0=D+k$, and $V_{\rm sh}=V+2R$. Then
$R\le V_{\rm sh}\le3R$, $A=K_0^2+V_{\rm sh}\ge R$, and $\Dom(A^{1/2})=H^1(\T)$.
The form inequality $K_0^2\le A$ gives \eqref{eq:Dres} and
$\|K_0A^{-1}\|^2\le R^{-1}$. Also
\[
 \|K_0^2A^{-1}\|=\|I-V_{\rm sh}A^{-1}\|\le1+3R\|A^{-1}\|\le4.
\]
For $f\in L^2(\T)$, $A^{-1}f\in H^2(\T)$ and
$K_0A^{-1}f\in H^1(\T)$, so the following form computation is legitimate:
\[
 \|A^{1/2}K_0A^{-1}f\|^2
 \le\|K_0^2A^{-1}f\|^2+3R\|K_0A^{-1}f\|^2
 \le19\|f\|^2.
\]
No derivative of $V$ is used.
\end{proof}

\subsection{Period-uniform form of the Floquet bound}

\begin{lemma}[{cf.\ \cite[Lemma 3.6]{Young}}]\label{lem:M2}
Let $s>2$ and $\psi\in\mathcal D_s$. There is $C_{\psi,s}$, independent of the period,
such that for every $p\ge\pi$, every $p$-periodic $V$ with $\|V\|_\infty\le R$,
and a.e.\ $k$,
\[
\big\|(\Hk+2R)U\psi(k,\cdot)\big\|^2_{L^2(\T,\dd q)}\le(1+3R)^2C_{\psi,s} .
\]
\end{lemma}

\begin{proof}
Since
$\|(\Hk+2R)\varphi\|\le(1+3R)(\|\Hz(k)\varphi\|+\|\varphi\|)$, it suffices to
bound $\|\Hz(k)U\psi(k,\cdot)\|$ and $\|U\psi(k,\cdot)\|$ uniformly. By
\eqref{eq:direct} with $V=0$, applied to $\psi\in H^2$, for a.e.\ $k$
\[
\Hz(k)U\psi(k,q)=-\sum_{a\in\Z}e^{-ik(q+pa)}\psi''(q+pa)\quad\text{in }L^2(\T,\dd q).
\]
The series converges absolutely for a.e.\ $q$, since $\psi''\in\Dom(|X|^s)$ with
$s>2$ \cite[Lem.\ 3.5]{Young}. After discarding the phases, monotone convergence and the
Cauchy--Schwarz inequality on each block give
\[
\begin{aligned}
\int_\T\big|\Hz(k)U\psi\big|^2\dd q
&\le\sum_{a,b\in\Z}\int_0^p|\psi''(q+pa)||\psi''(q+pb)|\,\dd q\\
&\le\Bigg(\sum_{a\in\Z}\Big(\int_{pa}^{p(a+1)}|\psi''|^2\Big)^{1/2}\Bigg)^2 .
\end{aligned}
\]
Choose an integer $N=N(\psi'',s)\ge1$ such that $\int_{|t|\ge x}(1+|t|^{2s})|\psi''|^2\le1$ for
$x\ge N$; it depends only on $\psi$ and $s$. For a block $I_a=[pa,p(a+1)]$ with
$|a|>N+1$, every $t\in I_a$ obeys
\[
|t|\ \ge\ p(|a|-1)\ \ge\ |a|-1\ \ge\ N,
\]
because $p\ge\pi>1$; note that the shift by one is needed for $a<0$, where the
block nearest the origin has $|t|\ge p(|a|-1)$ rather than $p|a|$. Hence
\[
\int_{I_a}|\psi''|^2\le\int_{|t|\ge p(|a|-1)}|\psi''|^2
\le\frac{1}{1+p^{2s}(|a|-1)^{2s}}\le\frac{1}{1+(|a|-1)^{2s}} .
\]
Bounding the finitely many blocks with $|a|\le N+1$ by $\|\psi''\|_{L^2(\R)}$
each and using $(1+y^2)^{-1/2}\le\sqrt2(1+y)^{-1}$ for the tail,
\[
\begin{aligned}
&\sum_{a\in\Z}\Big(\int_{pa}^{p(a+1)}|\psi''|^2\Big)^{1/2}\\
&\quad\le(2N+3)\|\psi''\|_{L^2(\R)}+\sum_{|a|>N+1}\frac{\sqrt2}{1+(|a|-1)^{s}}\\
&\quad\eqdef S_1(\psi,s)<\infty .
\end{aligned}
\]
Here $N$, $\|\psi''\|_{L^2}$ and the tail sum are independent of $p$; the period
entered only through the factor $p^{2s}\ge1$, which was discarded. The same
argument with $\psi$ in place of $\psi''$ defines a finite constant $S_0(\psi,s)$,
independent of $p$, such that $\|U\psi(k,\cdot)\|\le S_0(\psi,s)$. Take
$C_{\psi,s}\eqdef\big(S_0(\psi,s)+S_1(\psi,s)\big)^2$.
\end{proof}

\subsection{The master inequality}

\begin{proposition}\label{prop:master}
Let $R>0$, $s>2$, and $p\ge\pi$. Let $H$ be $p$-periodic with
$\|V\|_\infty\le R$, let $\psi\in\mathcal D_s$, and let
$C_{\psi,s}$ be as in Lemma~\ref{lem:M2}. For every $t,\e>0$,
\begin{equation}\label{eq:master}
 \left\|\mathcal A_H(t)\psi-\tfrac12Q_H\psi\right\|^2
 \le \frac{(1+3R)^2C_{\psi,s}}{R}
 \left[\frac{2p}{\pi}\mes(B_\e)+\frac{4\pi^2M_1}{3t^2\e^2}\right],
 \quad M_1=19.
\end{equation}
Here $B_\e$ is the bad set in Proposition~\ref{prop:badY}.
\end{proposition}

\begin{proof}
Write $u(k)=U\psi(k,\cdot)$, $A(k)=H(k)+2R$, and $K_0(k)=D+k$.
All fiber calculations below are made off the null set of symmetry points,
where the projections $P_m(k)$ are rank one. Their local analytic
representatives are measurable; arbitrary endpoint values lie on a null set.

\emph{Step 1: direct-integral and spectral representations.}
By \eqref{eq:direct}, \eqref{eq:intertwine}, and functional calculus,
\[
 UD(r)U^*=\int_{\T^*}^{\oplus}
 e^{irH(k)}K_0(k)e^{-irH(k)}\,\dk.
\]
The domain assertion follows from \eqref{eq:Dres}:
\[
 \int_{\T^*}\|K_0(k)e^{-irH(k)}u(k)\|^2\,\dk
 \le\|(H+2R)^{1/2}\psi\|^2<\infty.
\]
Put $\mathcal V(k)=\sum_mP_m(k)K_0(k)P_m(k)$. If $v_m$ is a unit eigenvector and
$\mu_m=E_m(k)+2R$, then
\[
 |\langle v_m,K_0(k)v_m\rangle|^2
 \le\|K_0(k)v_m\|^2\le\mu_m.
\]
Consequently, for $u\in\Dom(A^{1/2})$,
\begin{equation}\label{eq:formdomination}
 \|\mathcal V(k)u\|\le\|A^{1/2}u\|,\qquad
 \|e^{irH(k)}K_0(k)e^{-irH(k)}u\|\le\|A^{1/2}u\|.
\end{equation}
By \eqref{eq:GV}, $U(\tfrac12Q_H\psi)(k)=\mathcal V(k)u(k)$.

For fixed $k$, suppress $k$ in the notation and set
\[
 f_{mj}\eqdef\frac1t\int_0^t e^{i(E_m-E_j)r}\,\dd r.
\]
After the diagonal terms cancel, the fiber of
$\mathcal A_H(t)\psi-\tfrac12Q_H\psi$ is
\begin{equation}\label{eq:Y}
 Y(k)=\sum_m P_m\sum_{j\ne m}f_{mj}K_0P_ju.
\end{equation}
We justify both spectral sums. If $P_{>N}=I-\sum_{j\le N}P_j$, then
\[
 \sup_{0\le r\le t}\|K_0e^{-irH(k)}P_{>N}u\|
 \le\|A^{1/2}P_{>N}u\|\longrightarrow0.
\]
This proves convergence of the input sum uniformly in $r$; it does not assert
commutation of $K_0$ with $P_j$. For the output sum, its partial sums are
$P_{\le M}(e^{irH(k)}K_0e^{-irH(k)}-\mathcal V(k))u$, dominated in norm by
$2\|A^{1/2}u\|$. This bound is also square-integrable in $k$.
Dominated convergence for Bochner integrals \cite{HP} therefore justifies the
time integration and passage to the direct integral in \eqref{eq:Y}.
In particular, $k\mapsto\|Y(k)\|^2$ is measurable.

\emph{Step 2: a row estimate.}
Let $w=Au$ and $\mathcal T=A^{1/2}K_0A^{-1}$, so $\|\mathcal T\|^2\le M_1$ by
Proposition~\ref{prop:M1}. For
$S_m=\sum_{j\ne m}f_{mj}K_0P_ju$, the vector $S_m$ belongs to
$\Dom(A^{1/2})$ and
\[
 A^{1/2}S_m=\mathcal T\sum_{j\ne m}f_{mj}P_jw.
\]
Orthogonality, followed by Cauchy--Schwarz and the bound on $\mathcal T$, gives
\begin{align*}
 \|Y(k)\|^2
 &=\sum_m|\langle S_m,v_m\rangle|^2\\
 &\le\sum_m\mu_m^{-1}\|A^{1/2}S_m\|^2\\
 &\le M_1\sum_m\mu_m^{-1}
             \sum_{j\ne m}|f_{mj}|^2\|P_jw\|^2.
\end{align*}
In particular,
\begin{equation}\label{eq:316}
 \|Y(k)\|^2\le M_1\sum_m\frac1{\mu_m}
             \sum_{j\ne m}|f_{mj}|^2\|P_jw\|^2.
\end{equation}
The computation can first be made with finite output sums and then passed to
the limit; the right-hand side consists of nonnegative terms.

\emph{Step 3: the bad set.}
Equation~\eqref{eq:formdomination} directly yields
\begin{equation}\label{eq:badpointwise}
 \|Y(k)\|^2\le4\|A^{1/2}u\|^2\le\frac4R\|w\|^2.
\end{equation}
This estimate holds for almost every $k$, with no sum of resolvent weights.

\emph{Step 4: the good set.}
For $k\notin B_\e$, all adjacent gaps exceed $\e$. Thus
$|E_m-E_j|\ge|m-j|\e$ and
$|f_{mj}|\le2/(t\e|m-j|)$ for $m\ne j$. Since $\mu_m\ge R$, Tonelli's
theorem and \eqref{eq:316} imply
\begin{align}
 \|Y(k)\|^2
 &\le\frac{4M_1}{Rt^2\e^2}
     \sum_j\|P_jw\|^2\sum_{m\ne j}|m-j|^{-2}\notag\\
 &\le\frac{4\pi^2M_1}{3Rt^2\e^2}\|w\|^2.
 \label{eq:goodpointwise}
\end{align}
Here $\sum_{m\ge1,m\ne j}|m-j|^{-2}\le2\sum_{a\ge1}a^{-2}=\pi^2/3$.

\emph{Step 5: integration.}
Lemma~\ref{lem:M2} gives $\|w\|^2\le(1+3R)^2C_{\psi,s}$ uniformly in $k$.
Integrating \eqref{eq:badpointwise} over $B_\e$ and
\eqref{eq:goodpointwise} over its complement, with
$\dk=(p/2\pi)\dd k$ and total mass one, proves \eqref{eq:master}.
\end{proof}

\subsection{One-scale and two-scale consequences}
Throughout this subsection, $R>0$, $s>2$, $\psi\in\mathcal D_s$, and the
approximants satisfy $\|V_n\|_\infty\le R$ and
$p_{n+1}/p_n\in\N\setminus\{1\}$. The same $C_{\psi,s}$ is used at
all periods. Set
\begin{equation}\label{eq:Kprime}
 K'=K'(R)\eqdef\frac{(1+3R)^2}{R}
 \max\left\{\frac{10\sqrt{C_2\pi}}{\pi},\frac{4\pi^2M_1}{3}\right\}.
\end{equation}

\begin{corollary}\label{cor:EI}
For $p_n\ge\pi$ and $t,\e>0$,
\begin{equation}\label{eq:EI}
 \left\|\mathcal A_{H_n}(t)\psi-\tfrac12Q_n\psi\right\|^2
 \le K'C_{\psi,s}\left[p_ne^{cp_n}\sqrt\e+\frac1{t^2\e^2}\right].
\end{equation}
\end{corollary}
\begin{proof}
Insert $\mes(B_\e)\le4\sqrt{C_2\pi}e^{cp_n}\sqrt\e$ from
Proposition~\ref{prop:badY} into \eqref{eq:master} at period $p_n$.
The coefficient of $p_ne^{cp_n}\sqrt\e$ is bounded by the choice of $K'$.
\end{proof}

\begin{corollary}\label{cor:EII}
For $p_n\ge\pi$, $t>0$, and $\e>\gamma_{n+1}$,
\begin{equation}\label{eq:EII}
 \left\|\mathcal A_{H_{n+1}}(t)\psi-\tfrac12Q_{n+1}\psi\right\|^2
 \le K'C_{\psi,s}\left[p_{n+1}e^{cp_n}\sqrt\e+\frac1{t^2\e^2}\right].
\end{equation}
\end{corollary}
\begin{proof}
Apply \eqref{eq:master} at period $p_{n+1}$ and use
Theorem~\ref{thm:badset}. The coefficient of the bad-set term is
$10\sqrt{C_2\pi}/\pi$, as included in \eqref{eq:Kprime}.
\end{proof}

\begin{corollary}\label{cor:heis}
For $n\ge2$, $p_{n-1}\ge\pi$, $t>0$, and $\e>\gamma_n$,
\begin{equation}\label{eq:heis}
 \left\|Q_n\psi-t^{-1}X_{H_n}(t)\psi\right\|
 \le\frac{\|X\psi\|}{t}
 +2\left(K'C_{\psi,s}\left[p_ne^{cp_{n-1}}\sqrt\e+
           \frac1{t^2\e^2}\right]\right)^{1/2}.
\end{equation}
\end{corollary}
\begin{proof}
Divide \eqref{eq:33} by $t$ and use Corollary~\ref{cor:EII} with indices
shifted down by one.
\end{proof}

\begin{remark}[A balanced one-period rate]\label{rem:balanced}
For a single $p$-periodic operator, put
$\e=(pe^{cp})^{-2/5}t^{-4/5}$ in \eqref{eq:EI}. The two contributions are
equal, and hence
\[
 \left\|\mathcal A_H(t)\psi-\tfrac12Q_H\psi\right\|^2
 \le2K'C_{\psi,s}(pe^{cp})^{4/5}t^{-2/5}.
\]
Its square root has period factor $p^{2/5}e^{2cp/5}$ and time factor $t^{-1/5}$.
This improves the polynomial period factor in the rate of
\cite[Theorem 3.1]{Young}; its discriminant constant is kept distinct
from the explicit $c$ used here.
\end{remark}

\section{Proof of Theorem~\ref{thm:main}}\label{sec:main}
Throughout this section, $\eta>14c$, $H\in\EC(\eta,2)$,
$\|V\|_\infty\le R$, $s>2$, and $\psi\in\mathcal D_s$ are fixed.
The approximants are normalized as in Section~\ref{sec:fold}.
Write $C_{\psi,s}$ and $C_{\rm lin}=C_{\rm lin}(R,\psi)$ for the constants in
Lemma~\ref{lem:M2} and Proposition~\ref{prop:prop}. Theorem~\ref{thm:AK}
applies to each $H_n$, since $\mathcal D_s\subset H^1(\R)\cap\Dom(X)$.
The convergence, identification, and nonvanishing scheme follows
\cite[Section 4, proof of Theorem 1.2]{Young}, with the two-scale comparison
schedule of \cite[Sections 4--5]{HGS}. The estimates below use the periodic
rate and linear dynamical comparison established here.

\subsection{The dynamical inputs}
The geometry has supplied the periodic rate. Two additional facts are
needed to pass to the limit: a comparison linear in the potential error,
and a positive amount of Floquet mass away from symmetry momenta. Their
proofs are deferred to Appendix~\ref{app:dynamics}.
\begin{proposition}[Linear dynamical comparison]\label{prop:prop}
Let $H_i=-\partial_x^2+V_i$ with real $V_i\in L^\infty(\R)$ and
$\|V_i\|_\infty\le R$. For $\psi\in H^2(\R)\cap\Dom(X)$ and $t\in\R$,
\begin{equation}\label{eq:prop}
 \|X_{H_1}(t)\psi-X_{H_2}(t)\psi\|
 \le C_{\rm lin}(R,\psi)(|t|^{3/2}+t^2)\|V_1-V_2\|_\infty.
\end{equation}
One may take, with $m=\|\psi\|$,
\[
 h=\|\psi''\|+2Rm,\qquad F=(\|\psi'\|^2+3Rm^2)^{1/2},\qquad
 C_{\rm lin}=\frac{4\sqrt2}{3}m+\sqrt{mh}+F.
\]
The more precise bound is
$\|V_1-V_2\|_\infty[(4\sqrt2/3)m|t|^{3/2}+(\sqrt{mh}+F)t^2]$.
\end{proposition}
The proof is in Appendix~\ref{app:linear}. It uses the common $H^2$ domain
and does not differentiate either potential. Young's square-root comparison \cite[Lemma 4.2]{Young} suffices for
his next-period approximation class. The linear dependence here permits
the second-power weight in Theorem~\ref{thm:main}.

\begin{proposition}[{\cite[Lemma 4.3]{Young}}]\label{prop:mass}
Let $\psi\in H^1(\R)\cap\Dom(|X|^s)$ with $s>2$ and
\begin{equation}
    S_n\eqdef \pm\left[\frac{\pi}{4p_n},\frac{3\pi}{4p_n}\right].
\end{equation}
For the fixed period sequence $(p_n)$, there exists $N_0(\psi)$, also
depending on that sequence, such that
\begin{equation}
    \int_{S_n}\|U_{p_n}\psi(k,\cdot)\|^2\frac{\dd k}{|\T_n^\ast|}\geq\frac18\|\psi\|^2,
    \qquad n\geq N_0(\psi).
\end{equation}
\end{proposition}

\begin{remark}\label{rem:halfpower}
The common-domain identity in Appendix~\ref{app:linear} controls both the
state difference and its graph norm linearly in the potential error.
Using instead the square-root comparison of \cite[Lemma 4.2]{Young}
in the same two-scale argument gives the sufficient condition
$H\in\EC(\eta,4)$ with $\eta>28c$.
\end{remark}

\subsection{Choice of parameters}
Choose
\begin{equation}\label{eq:window}
 10c<\eta_0<\frac45(\eta-\tfrac32c),
\end{equation}
and put
\begin{align}
 \e_n&\eqdef p_n^{-2}e^{-\eta_0p_n},&
 \widetilde\e_n&\eqdef p_{n+1}^{-2}e^{-\eta_0p_n},\notag\\
 t_n&\eqdef p_{n+1}^2e^{\zeta p_n},&
 \zeta&\eqdef\frac{5\eta_0}{4}-\frac c2.
 \label{eq:choices}
\end{align}
Also define
\begin{equation}\label{eq:deltas}
 \delta\eqdef\frac{\eta_0}{4}-\frac c2,\qquad
 \delta'\eqdef\eta-\frac{5\eta_0}{4}+\frac c2.
\end{equation}

\begin{lemma}\label{lem:param}
The interval \eqref{eq:window} is nonempty precisely when $\eta>14c$.
The times $t_n$ increase to infinity and are eventually at least one.
For all sufficiently large $n$, $\gamma_{n+1}<\widetilde\e_n$ and
$\gamma_n<\widetilde\e_{n-1}$. Moreover, $\delta,\delta'>2c$, and
\begin{equation}\label{eq:balance}
 p_{n+1}e^{cp_n}\sqrt{\widetilde\e_n}
 =\frac1{t_n^2\widetilde\e_n^2}=e^{-2\delta p_n}.
\end{equation}
\end{lemma}
\begin{proof}
The inequality $10c<\tfrac45(\eta-\tfrac32c)$ is equivalent to $\eta>14c$.
Since $\zeta>0$ and the periods increase, $t_n\uparrow\infty$.
The estimate \eqref{eq:gamma} with $\al=2$ gives
\[
 \frac{\gamma_{n+1}}{\widetilde\e_n}
 \le2e^{-(\eta-\eta_0)p_n}\longrightarrow0,
\]
because $\eta_0<\eta$. The shifted admissibility assertion is identical.
The inequalities $\delta>2c$ and $\delta'>2c$ are respectively equivalent
to the lower and upper bounds in \eqref{eq:window}. Finally,
\[
 p_{n+1}e^{cp_n}\sqrt{\widetilde\e_n}
 =e^{(c-\eta_0/2)p_n},\qquad
 \frac1{t_n^2\widetilde\e_n^2}
 =e^{(-2\zeta+2\eta_0)p_n}=e^{(c-\eta_0/2)p_n},
\]
which proves \eqref{eq:balance}.
\end{proof}

The period power $2$ in $t_n$ cancels the second-power approximation weight
under the linear perturbation estimate. The strict bounds $\delta,\delta'>2c$ will ensure
that the limit of the periodic velocities is nonzero.

\subsection{Convergence of the periodic velocities}
\begin{proposition}\label{prop:cauchy}
For all sufficiently large $n$,
\begin{equation}\label{eq:cauchy}
 \|Q_n\psi-Q_{n+1}\psi\|
 \le4(2K'C_{\psi,s})^{1/2}e^{-\delta p_n}
       +4C_{\rm lin} e^{-\delta'p_n}.
\end{equation}
In particular, $\sum_n\|Q_n\psi-Q_{n+1}\psi\|<\infty$, and
$Q_\infty\psi\eqdef\lim_n Q_n\psi$ exists in $L^2(\R)$.
\end{proposition}
\begin{proof}
Use the same time $t_n$ for both approximants and write
\[
 \tfrac12\|Q_n\psi-Q_{n+1}\psi\|\le I_n+II_n+III_n,
\]
where
\begin{align*}
 I_n&=\|\mathcal A_{H_n}(t_n)\psi-\tfrac12Q_n\psi\|,\\
 II_n&=\|\mathcal A_{H_{n+1}}(t_n)\psi-\tfrac12Q_{n+1}\psi\|,\\
 III_n&=\|\mathcal A_{H_n}(t_n)\psi-\mathcal A_{H_{n+1}}(t_n)\psi\|.
\end{align*}
For $I_n$, Corollary~\ref{cor:EI} at $\e=\e_n$ gives
\[
 p_ne^{cp_n}\sqrt{\e_n}=e^{-2\delta p_n},\qquad
 \frac1{t_n^2\e_n^2}
 =\left(\frac{p_n}{p_{n+1}}\right)^4e^{-2\delta p_n}
 \le e^{-2\delta p_n}.
\]
Thus $I_n\le(2K'C_{\psi,s})^{1/2}e^{-\delta p_n}$.
For $II_n$, use Corollary~\ref{cor:EII} with
$\e=\widetilde\e_n$, admissible by Lemma~\ref{lem:param}.
Both contributions equal $e^{-2\delta p_n}$ by \eqref{eq:balance}, hence
$II_n\le(2K'C_{\psi,s})^{1/2}e^{-\delta p_n}$.

For the dynamical comparison, \eqref{eq:33} implies
\[
 III_n=\frac1{2t_n}\|X_{H_n}(t_n)\psi-X_{H_{n+1}}(t_n)\psi\|.
\]
Since $t_n\ge1$, Proposition~\ref{prop:prop} and \eqref{eq:gamma} give
\begin{align*}
 III_n
 &\le\frac{C_{\rm lin}}{2t_n}(t_n^{3/2}+t_n^2)\gamma_{n+1}
 \le C_{\rm lin} t_n\gamma_{n+1}\\
 &\le2C_{\rm lin} p_{n+1}^{2}p_{n+1}^{-2}e^{(\zeta-\eta)p_n}
 =2C_{\rm lin} e^{-\delta'p_n}.
\end{align*}
Adding the three bounds and multiplying by two proves \eqref{eq:cauchy}.
For any fixed sufficiently large $N$ and $n\ge N$, we have
$p_n\ge2^{n-N}p_N$. Since $\delta,\delta'>0$, the series converges.
\end{proof}

\begin{corollary}\label{cor:tail}
For all sufficiently large $n$, with
$C_Q\eqdef8(2K'C_{\psi,s})^{1/2}+8C_{\rm lin}$,
\begin{equation}\label{eq:tail}
 \|Q_n\psi-Q_\infty\psi\|\le C_Q(e^{-\delta p_n}+e^{-\delta'p_n}).
\end{equation}
\end{corollary}
\begin{proof}
For $k\ge n$, $p_k\ge2^{k-n}p_n\ge(k-n+1)p_n$, so
\[
 \sum_{k\ge n}e^{-\delta p_k}
 \le\frac{e^{-\delta p_n}}{1-e^{-\delta p_n}}
 \le2e^{-\delta p_n}
\]
for large $n$. Apply this to both exponents and sum \eqref{eq:cauchy}.
\end{proof}

\subsection{Identification for all large times}
\begin{proposition}\label{prop:limit}
One has $\lim_{t\to\infty}t^{-1}X_H(t)\psi=Q_\infty\psi$ in $L^2(\R)$.
\end{proposition}
\begin{proof}
For sufficiently large $t$, choose $n$ such that $t_{n-1}\le t<t_n$ and split
\begin{align*}
 \|t^{-1}X_H(t)\psi-Q_\infty\psi\|
 \le{}&t^{-1}\|X_H(t)\psi-X_{H_n}(t)\psi\|\\
 &+\|t^{-1}X_{H_n}(t)\psi-Q_n\psi\|+\|Q_n\psi-Q_\infty\psi\|.
\end{align*}
By Definition~\ref{def:class}, eventually
$\|V-V_n\|_\infty\le p_{n+1}^{-2}e^{-\eta p_n}$.
For the first term, Proposition~\ref{prop:prop} and $1\le t<t_n$ imply
\[
 t^{-1}\|X_H(t)\psi-X_{H_n}(t)\psi\|
 \le2C_{\rm lin} t_n\|V-V_n\|_\infty
 \le2C_{\rm lin} e^{-\delta'p_n}\longrightarrow0.
\]
For the second term, use Corollary~\ref{cor:heis} with
$\e=\widetilde\e_{n-1}$ and $t\ge t_{n-1}$. By the shifted form of
\eqref{eq:balance}, this term is at most
\[
 t_{n-1}^{-1}\|X\psi\|+2(2K'C_{\psi,s})^{1/2}e^{-\delta p_{n-1}},
\]
which tends to zero. The third term tends to zero by \eqref{eq:tail}.
\end{proof}
This identifies $Q_\infty\psi$ with $Q\psi$ in Theorem~\ref{thm:main}.
The limit is therefore independent of the chosen admissible sequence of
periodic approximants. It defines a linear map on $\bigcup_{s>2}\mathcal D_s$;
no self-adjoint extension of this map is needed below.

\subsection{Nonvanishing}
\begin{proposition}\label{prop:nonzero}
For every $\psi\in\mathcal D_s\setminus\{0\}$, one has $Q\psi\ne0$.
\end{proposition}
\begin{proof}
Let $S_n=\pm[\pi/(4p_n),3\pi/(4p_n)]$. Theorem~\ref{thm:AK}, unitarity of
$U_{p_n}$, and orthogonality of the rank-one projections give
\[
 \|Q_n\psi\|^2
 =\int_{\T_n^*}\sum_m |(E_m^{(n)})'(k)|^2
       \|P_m^{(n)}(k)U_{p_n}\psi(k)\|^2\,\dkn{n}.
\]
On $S_n$, $|\sin(p_nk)|\ge\sqrt2/2$, so \eqref{eq:26} gives
$|(E_m^{(n)})'(k)|\ge\sqrt2(C_2p_n)^{-1}e^{-2cp_n}$.
Consequently, Proposition~\ref{prop:mass} implies, for all large $n$,
\begin{equation}\label{eq:415}
 \|Q_n\psi\|^2
 \ge\frac{2e^{-4cp_n}}{C_2^2p_n^2}
       \int_{S_n}\|U_{p_n}\psi(k)\|^2\,\dkn{n}
 \ge\frac{e^{-4cp_n}}{4C_2^2p_n^2}\|\psi\|^2.
\end{equation}
Combining this with \eqref{eq:tail} yields
\[
 \|Q\psi\|\ge\frac{e^{-2cp_n}}{2C_2p_n}\|\psi\|
            -C_Q(e^{-\delta p_n}+e^{-\delta'p_n}).
\]
Because $\delta,\delta'>2c$, the subtracted terms are
$o(p_n^{-1}e^{-2cp_n})$. The right-hand side is therefore strictly positive
for some sufficiently large $n$.
\end{proof}

\begin{proof}[Proof of Theorem~\ref{thm:main}]
Propositions~\ref{prop:cauchy}, \ref{prop:limit}, and \ref{prop:nonzero}
establish the asserted strong limit and its nonvanishing.
\end{proof}

\subsection{All moment orders between zero and four}
\begin{proof}[Proof of Corollary~\ref{cor:moments}]
By Theorem~\ref{thm:main},
\begin{equation}\label{eq:secondmomentlimit}
 t^{-2}M_{2,\psi}(t)=\|t^{-1}X_H(t)\psi\|^2
 \longrightarrow\|Q\psi\|^2>0.
\end{equation}
Normalize $\|\psi\|=1$; rescaling the initial state changes only the
final constants. Equation~\eqref{eq:secondmomentlimit}
gives constants $a_*,A_*>0$ such that $a_*t^2\le M_{2,\psi}(t)\le A_*t^2$ for all large $t$.
Since $\psi\in H^2(\R)\cap\Dom(X^2)$, the fourth-moment estimate of
Lemma~\ref{lem:Bmoment} (see also \cite[Theorem 5.2]{Young}) gives $M_{4,\psi}(t)\le b_* t^4$ for large $t$, for some
$b_*>0$. In particular, all moments of order at most four are finite.

The measure $|e^{-itH}\psi(x)|^2\dd x$ is a probability measure. For
$2\le r\le4$, the standard comparison of moments gives
\[
 M_{2,\psi}(t)^{r/2}\le M_{r,\psi}(t)\le M_{4,\psi}(t)^{r/4},
\]
which proves \eqref{eq:ballisticmoments} in that range.
For $0<r<2$, H\"older's inequality with $\theta=2/(4-r)\in(0,1)$ yields
\[
 M_{2,\psi}(t)\le M_{r,\psi}(t)^{\theta}M_{4,\psi}(t)^{1-\theta}.
\]
Therefore
\begin{equation}\label{eq:lowermoments}
 M_{r,\psi}(t)\ge
 M_{2,\psi}(t)^{(4-r)/2}M_{4,\psi}(t)^{-(2-r)/2}
 \ge a_*^{(4-r)/2}b_*^{-(2-r)/2}t^r.
\end{equation}
The upper bound is $M_{r,\psi}(t)\le M_{2,\psi}(t)^{r/2}\le A_*^{r/2}t^r$.
It follows that
\eqref{eq:ballisticmoments} holds for every $0<r\le4$; dividing its logarithm
by $r\log t$ proves the asserted transport exponents.
\end{proof}

\section{Remarks and open questions}\label{sec:remarks}
\begin{remark}[The state class and higher moments]\label{ex:highmoments}
The endpoint $4$ in Corollary~\ref{cor:moments} concerns the regularity of
the initial state. For instance, let $\chi\in C_c^\infty(\R)$ equal one
near zero and put $\psi(x)=\chi(x)|x|^{9/4}$.
Then $\psi\in H^2\setminus H^3$ and $\psi\in\mathcal D_s$ for every
$s>0$. For the free operator and $t\ne0$,
\begin{equation}\label{eq:freecounterexample}
 (e^{-itH_0}\psi)(x)=c_t e^{ix^2/(4t)}
 \widehat{e^{iy^2/(4t)}\psi(y)}\bigl(x/(2t)\bigr),\qquad c_t\ne0.
\end{equation}
Plancherel shows that a finite sixth position moment would require
$e^{iy^2/(4t)}\psi\in H^3$. The phase and its inverse preserve Sobolev
regularity on the compact support, so $M_{6,\psi}(t)=\infty$.

More generally, for any $r>4$ choose
\[
 \beta=\frac32+\min\{(r-4)/4,1/8\},\qquad
 \psi_r=\chi|x|^\beta,\qquad \theta_0=\beta-\frac32\in(0,1).
\]
Then $\psi_r\in H^2$ and $\beta+1/2<r/2$.
The Fourier--Gagliardo characterization \cite[Proposition 3.4]{DPV}
shows that $\psi_r''\notin H^{\theta_0}$: on a sufficiently small
positive interval its seminorm contains, up to a positive factor,
\[
 \int_0^\varepsilon\!\int_0^x
 \frac{|x^{\beta-2}-y^{\beta-2}|^2}{|x-y|^{1+2\theta_0}}\,\dd y\,\dd x
 =\left(\int_0^1
 \frac{|1-u^{\beta-2}|^2}{(1-u)^{1+2\theta_0}}\,\dd u\right)
 \int_0^\varepsilon\frac{\dd x}{x}=\infty.
\]
The parenthesized integral is positive and finite because
$-1/2<\beta-2<0$ and $0<\theta_0<1$.
Thus $\psi_r\notin H^{2+\theta_0}$ and consequently
$\psi_r\notin H^{r/2}$, whereas compact support puts it in every
$\mathcal D_s$. Formula~\eqref{eq:freecounterexample} now gives
$M_{r,\psi_r}(t)=\infty$ for every $t\ne0$.
\end{remark}

\begin{remark}[Extension to the energy space]\label{bd:h1extension}
The limit in Theorem~\ref{thm:main} extends uniquely to a bounded linear
map $Q_{\rm en}:H^1(\R)\to L^2(\R)$. Indeed, form invariance gives
\[
 \sup_{t>0}\|2\mathcal A_H(t)\psi\|
 \le2\|(H+2R)^{1/2}\psi\|\le C_R\|\psi\|_{H^1}.
\]
On $\bigcup_{s>2}\mathcal D_s$ the same bound holds for $Q$ by passage
to the limit. This space contains $C_c^\infty(\R)$ and is dense in
$H^1$. Approximation in $H^1$ therefore gives
$2\mathcal A_H(t)\psi\to Q_{\rm en}\psi$ for every $\psi\in H^1$.
The Heisenberg identity then gives
$t^{-1}X_H(t)\psi\to Q_{\rm en}\psi$ on $H^1\cap\Dom(X)$.
Injectivity is established on the original weighted class only; the
nonvanishing estimate does not provide a lower bound on the $H^1$ unit
sphere.
\end{remark}

\begin{remark}[Open quantitative questions]\label{rem:a9}
The approximation power $2$ and threshold $14c$ are sufficient.
The present choice of parameters makes the periodic errors summable and
ensures $t_n\gamma_{n+1}=o(p_n^{-1}e^{-2cp_n})$, the scale needed for
nonvanishing. Optimal approximation conditions remain to be determined;
see \cite{DF2,HGS}.
One can also compare $H_{n+1}$ with $H_{n-j}$ in the large cell, but the
relevant perturbation is
$\sum_{i=0}^j\gamma_{n+1-i}$; the next-scale envelope alone need not make
this sum smaller than the chosen separation parameter. A useful deeper
comparison requires additional error control or different parameters.
Finally, higher moments require stronger momentum regularity together
with propagation of the corresponding position domains, as the preceding
free examples demonstrate.
\end{remark}

\appendix
\section{The uniform discriminant bound}\label{app:disc}

\begin{proof}[Proof of Lemma~\ref{lem:discriminant}]
Fix an energy $E\ge-R$. For a positive number $s_0$, kept fixed as the
energy parameter $\lambda$ varies, write the equation in the state
$(u,u'/s_0)^{\mathsf T}$:
\[
 Y'=B_\lambda(x)Y,\qquad
 B_\lambda(x)=\begin{pmatrix}0&s_0\\(V(x)-\lambda)/s_0&0\end{pmatrix},
 \quad
 \partial_\lambda B_\lambda =\begin{pmatrix}0&0\\-s_0^{-1}&0\end{pmatrix}.
\]
Let $\Phi_\lambda(x,y)$ be its fundamental matrix. The fixed conjugation
from $(u,u')$ does not change its monodromy trace. Differentiating the
Volterra equation gives
\[
 \Delta'(E)=\tr\int_0^p
 \Phi_E(p,x)\partial_\lambda B_\lambda\Phi_E(x,0)\,\dd x.
\]
The formula follows from uniformly convergent Volterra iterations on a
bounded cell; it requires only bounded measurable $V$. Whenever
$\|\Phi_E(x,y)\|\le e^{a_0(x-y)}$ for $x\ge y$, it implies
\begin{equation}\label{eq:frozen-first}
 |\Delta'(E)|\le2(p/s_0)e^{a_0p}.
\end{equation}
Only one exponential appears because the two propagation lengths add to
$p$.

For $-R\le E\le R$, choose
$s_0=\max\{p^{-1},\sqrt{2R}\}$. Since $|V-E|\le2R\le s_0^2$,
$\|B_E(x)\|\le s_0$. Gr\"onwall's inequality and \eqref{eq:frozen-first} give
\[
 |\Delta'(E)|\le2p^2e^{\max\{1,\sqrt{2R}\,p\}}
 \le2e\,p^2e^{\sqrt{2R}\,p}.
\]
For $E\ge R$ and $E\ge p^{-2}$, choose $s_0=\sqrt E$ and write
\[
 B_E=s_0\begin{pmatrix}0&1\\-1&0\end{pmatrix}
       +\frac{V}{s_0}\begin{pmatrix}0&0\\1&0\end{pmatrix}.
\]
Conjugating out the orthogonal rotation leaves norm at most $R/s_0$.
Thus
\[
 |\Delta'(E)|\le2(p/s_0)e^{Rp/s_0}
 \le2p^2e^{\sqrt R\,p}.
\]
Finally, if $R\le E<p^{-2}$, choose $s_0=p^{-1}$. Then
$|V-E|<2p^{-2}$, so $\|B_E(x)\|\le2/p$ and
$|\Delta'(E)|\le2e^2p^2$. The three cases cover $E\ge-R$ and are bounded
by $2e^3p^2e^{\sqrt{2R}\,p}$. Since $2c=\sqrt{2R}$, this is
\eqref{eq:disc}.
\end{proof}

\section{Dynamical and fiber-mass inputs}\label{app:dynamics}
\subsection*{The one-scale bad set}\label{app:badone}
\begin{proof}[Proof of Proposition~\ref{prop:badY}]
If $m<m'$ and $|E_m(k)-E_{m'}(k)|\le\varepsilon$, non-overlap gives
$E_{m'}(k)\le\max E_m+\varepsilon\le\min E_{m'}+\varepsilon$.
Choose the nearest minimizing symmetry point $k_*$, interpreted modulo the
ends of the Brillouin zone. Its distance $q$ from $k$ lies in $[0,\pi/p]$.
The one-period estimate \eqref{eq:27}, applied between $k$ and that
minimizing endpoint using evenness and periodicity, gives
$\varepsilon\ge e^{-2cp}q^2/(\pi C_2)$, and therefore
$q\le\rho=\sqrt{\pi C_2}e^{cp}\sqrt\varepsilon$.
Every colliding pair is thus contained in the same union of intervals
$[-\rho,\rho]$, $[\pi/p-\rho,\pi/p]$, and $[-\pi/p,-\pi/p+\rho]$,
intersected with the zone. Their total length is at most $4\rho$.
\end{proof}

\subsection*{The Heisenberg identity and fourth moment}
The domain and moment estimates originate in
\cite[Theorems 2.1--2.2]{RadSim}; their bounded-potential formulation is
\cite[Theorems 5.1--5.2]{Young}. We retain proofs to track the explicit
constants and our normalization.
Throughout, $V$ is real, $\|V\|_\infty\le R$, and $\phi_t=e^{-itH}\psi$.
\begin{lemma}[{cf.\ \cite[equation (3.3)]{Young}}]\label{lem:BHeis}
For $\psi\in H^1(\R)\cap\Dom(X)$, $\phi_t\in\Dom(X)$ for all $t$, and
\eqref{eq:33} holds. Moreover
\[
 \|X\phi_t\|\le\|X\psi\|+2|t|\|(H+2R)^{1/2}\psi\|.
\]
\end{lemma}
\begin{proof}
Take $F_\varepsilon(x)=x/(1+\varepsilon x^2)$.
It is bounded and smooth, $|F_\varepsilon'|\le1$,
$F_\varepsilon'\to1$ pointwise, $\|F_\varepsilon''\|_\infty\to0$, and
$|F_\varepsilon(x)|\uparrow|x|$ as $\varepsilon\downarrow0$.
For $\psi\in H^2$ the bounded-multiplier identity is
\[
 e^{itH}F_\varepsilon(X)e^{-itH}\psi
 =F_\varepsilon(X)\psi+
 \int_0^te^{irH}(2F_\varepsilon'D-iF_\varepsilon'')e^{-irH}\psi\,\dd r.
\]
For fixed $\varepsilon$ it extends by density to $H^1$, since the commutator
maps $H^1$ boundedly to $L^2$ and the evolution preserves the form domain.
The form inequality
\[
 \|D\phi_r\|^2\le\|(H+2R)^{1/2}\phi_r\|^2
 =\|(H+2R)^{1/2}\psi\|^2
\]
provides an integrable bound on compact time intervals. For an initial vector
also in $\Dom(X)$, the right side therefore converges in $L^2$ to
$Z=X\psi+2\int_0^tD(r)\psi\,\dd r$. Monotone convergence gives
$\int x^2|\phi_t|^2=\lim\|F_\varepsilon\phi_t\|^2=\|Z\|^2<\infty$.
The limiting vector is $e^{itH}X\phi_t$, proving the identity and the bound.
\end{proof}

\begin{lemma}[Fourth-moment estimate; cf.\ {\cite[Theorem 5.2]{Young}}]\label{lem:Bmoment}
If $\psi\in H^2(\R)\cap\Dom(X^2)$, then $\phi_t\in\Dom(X^2)$ and, with
$h_{\rm mom}=\|H\psi\|+R\|\psi\|$,
\[
 \|X^2\phi_t\|\le
 \left((1+\|X^2\psi\|^2)^{1/4}+(3\|\psi\|+2\sqrt{h_{\rm mom}})|t|\right)^2.
\]
In particular $M_{4,\psi}(t)\le b t^4$ for $|t|\ge1$, uniformly over
$\|V\|_\infty\le R$ when $\|\psi\|$, $\|\psi''\|$, and $\|X^2\psi\|$ are fixed.
\end{lemma}
\begin{proof}
Put $G_\varepsilon=x^2/(1+\varepsilon x^2)$,
$w_\varepsilon=(G_\varepsilon')^2$, and $g(t)=\|G_\varepsilon\phi_t\|^2$.
Then $0\le G_\varepsilon\uparrow x^2$, $|G_\varepsilon''|\le2$,
$w_\varepsilon\le4G_\varepsilon$, and
\[
 w_\varepsilon''=\frac{8(1-18y+21y^2)}{(1+y)^6},\qquad y=\varepsilon x^2.
\]
To verify $|w_\varepsilon''|\le8$, put $\mathfrak p(y)=1-18y+21y^2$.
For $y\ge0$,
\begin{align*}
 (1+y)^6-\mathfrak p(y)&=y(24-6y+20y^2+15y^3+6y^4+y^5)\ge0,\\
 (1+y)^6+\mathfrak p(y)&=36(y-1/6)^2+1+20y^3+15y^4+6y^5+y^6>0.
\end{align*}
In the first line $20y^2-6y+24>0$ because its discriminant is negative.
This proves the uniform multiplier bound.
The evolution preserves $H^2$, and
$\|\phi_t''\|\le\|H\psi\|+R\|\psi\|=h_{\rm mom}$.
The bounded commutator gives
$|g'|\le2\sqrt g(2\|G_\varepsilon'\phi_t'\|+2\|\psi\|)$.
Integration by parts, justified first for compactly supported smooth
approximants in $H^2$, gives
\[
 \|G_\varepsilon'\phi_t'\|^2
 =\frac12\int w_\varepsilon''|\phi_t|^2
       -\operatorname{Re}\int w_\varepsilon\overline{\phi_t}\phi_t''
 \le4\|\psi\|^2+4\sqrt g\,h_{\rm mom}.
\]
All multipliers and their derivatives are bounded for fixed $\varepsilon$,
so $H^2$ approximation passes to this identity. Hence
$\|G_\varepsilon'\phi_t'\|\le2\|\psi\|+2\sqrt{h_{\rm mom}}\,g^{1/4}$.
Set $Y=(1+g)^{1/4}\ge1$. Then
$2Y|Y'|\le6\|\psi\|+4\sqrt{h_{\rm mom}}Y$, and
$|Y'|\le3\|\psi\|+2\sqrt{h_{\rm mom}}$.
Integrate in either time direction and let $\varepsilon\downarrow0$.
Monotone convergence proves both the domain assertion and the bound.
\end{proof}

\subsection*{Linear comparison}\label{app:linear}
\begin{proof}[Proof of Proposition~\ref{prop:prop}]
Write $U_i(t)=e^{-itH_i}$, $\Delta V=V_1-V_2$, $d=\|\Delta V\|_\infty$,
$T=|t|$, and $u(t)=(U_1(t)-U_2(t))\psi$. Both operators have domain
$H^2(\R)$, which their own propagators preserve. The Duhamel estimate \cite[Lemma 4.1]{Young} gives
\[
 \|u(t)\|\le Td\,m,\qquad m=\|\psi\|.
\]
The decisive identity, with every term in $L^2$, is
\begin{align*}
 H_1u(t)
 &=U_1(t)H_1\psi-(H_2+\Delta V)U_2(t)\psi\\
 &=(U_1(t)-U_2(t))H_2\psi+U_1(t)(\Delta V)\psi-(\Delta V)U_2(t)\psi.
\end{align*}
This uses $(\Delta V)\psi$ only in $L^2$; no derivative or domain invariance of $\Delta V$
is assumed. Therefore
\[
 \|H_1u(t)\|\le d(T\|H_2\psi\|+2m).
\]
Using $H_1=D^2+V_1$ and $\|V_1\|_\infty\le R$,
\begin{align*}
 \|Du(t)\|^2
 &\le\|u(t)\|\,\|H_1u(t)\|+R\|u(t)\|^2\\
 &\le d^2\{2Tm^2+T^2m(\|H_2\psi\|+Rm)\}\\
 &\le d^2(2Tm^2+T^2mh),\qquad h=\|\psi''\|+2Rm.
\end{align*}
In particular, $\|Du(t)\|\le d(\sqrt{2T}\,m+T\sqrt{mh})$.
Both factors in the quadratic-form estimate are small in $d$; taking a square root does not lose linear dependence.

Next $D_i(r)=U_i(-r)DU_i(r)$, so
\[
 (D_1(r)-D_2(r))\psi
 =U_1(-r)D(U_1(r)-U_2(r))\psi
 +(U_1(-r)-U_2(-r))DU_2(r)\psi.
\]
The vector $DU_2(r)\psi$ belongs to $L^2$ and has norm at most
$F=(\|\psi'\|^2+3Rm^2)^{1/2}$. Duhamel applies to that vector too, yielding
\[
 \|(D_1(r)-D_2(r))\psi\|
 \le d\{\sqrt{2|r|}\,m+|r|(\sqrt{mh}+F)\}.
\]
Apply Lemma~\ref{lem:BHeis} to each operator and subtract. The initial terms
$X\psi$ cancel. Bochner integration gives
\begin{align*}
 \|X_{H_1}(t)\psi-X_{H_2}(t)\psi\|
 &\le2d\int_0^T\{\sqrt{2r}\,m+r(\sqrt{mh}+F)\}\,\dd r\\
 &=d\left\{\frac{4\sqrt2}{3}mT^{3/2}+(\sqrt{mh}+F)T^2\right\}.
\end{align*}
This proves \eqref{eq:prop}, also at $t=0$. The constants are uniform in
both potentials subject to the given bound. The stated domain $H^2\cap\Dom(X)$ is sufficient throughout.
\end{proof}

\subsection*{Fiber mass away from the symmetry points}
\begin{proof}[Proof of Proposition~\ref{prop:mass}]
The limiting assertion holds for every $\psi\in L^2(\R)$.
First let $\psi$ have compact support and choose $p$ larger than twice
its support radius. On $0\le q<p$, only the cells $0,-1$ contribute:
\[
 U_p\psi(k,q)=e^{-ikq}\{\psi(q)+e^{ikp}\psi(q-p)\}.
\]
For $S_p=\pm[\pi/(4p),3\pi/(4p)]$,
\[
 \int_{S_p}\frac{p\,\dd k}{2\pi}=\frac12,\qquad
 \int_{S_p}e^{ikp}\,\dd k
 =\frac2p\bigl(\sin(3\pi/4)-\sin(\pi/4)\bigr)=0.
\]
Expanding the square and integrating therefore gives
$\|\mathbf1_{S_p}U_p\psi\|^2=\tfrac12\|\psi\|^2$.
For general $\psi$, approximate by compactly supported $\varphi$.
The operators $\mathbf1_{S_p}U_p$ are contractions, so
\[
 \bigl|\|\mathbf1_{S_p}U_p\psi\|^2-
        \|\mathbf1_{S_p}U_p\varphi\|^2\bigr|
 \le\|\psi-\varphi\|(\|\psi\|+\|\varphi\|)
\]
uniformly in $p$. Letting $p\to\infty$ and then
$\|\psi-\varphi\|\to0$ proves convergence to $\tfrac12\|\psi\|^2$.
The stated lower bound follows for $p=p_n$ sufficiently large.
\end{proof}

\section*{Acknowledgments}
The author thanks Rui Han for helpful discussions, including discussions of
the literature on limit-periodic potentials, and Bart Rosenzweig for comments
on the manuscript and suggestions concerning the exposition.
Special thanks are due to Giorgio Young for his comments on an earlier
version of the manuscript, which led to improvements in the exposition.

\subsection*{Funding}
M. Gomez Solis was supported in part by the National Science Foundation under
grant DMS-2143369.

\subsection*{AI-assisted manuscript preparation}
OpenAI’s ChatGPT was used for the preceding mathematical review and for preparing this revision. The prompts requested a statement-by-statement
proof audit. Most notably, this audit led to the interpolation argument used in the proof of Collorary~\ref{cor:moments}. Responsibility
for the mathematics, citations, declarations, and final submitted text
remains with the author.

\subsection*{Competing interests}
The author has no relevant financial or non-financial interests to disclose.

\subsection*{Data availability}
The results are theoretical and are established by the mathematical
arguments in this article. No empirical dataset is required for these arguments.

\enlargethispage{3\baselineskip}

\end{document}